\documentclass[english,reqno]{amsart}

\usepackage[margin=1.15in]{geometry}
\usepackage{amsmath,amssymb,mathtools}
\usepackage{mathrsfs}
\usepackage{dsfont}
\usepackage{xcolor}
\usepackage{hyperref}
\usepackage[nameinlink]{cleveref}

\definecolor{linkblue}{rgb}{0,0.447,0.698}

\usepackage{pgfplots}
\pgfplotsset{compat=1.18}
\usepgfplotslibrary{fillbetween}

\hypersetup{
  colorlinks=true,
  linkcolor=linkblue,
  citecolor=linkblue,
  urlcolor=linkblue,
  filecolor=linkblue,
  breaklinks=true,
  bookmarksnumbered=true,
  pdftitle={Sharp boundary regularity for planar minimal graphs with analytic Dirichlet data},
  pdfauthor={Camillo De Lellis, Guido De Philippis, Alessio Figalli},
  pdfsubject={},
  pdfkeywords={Minimal graph, least area problem, obstacle problem, free boundary, boundary regularity}}

  \newcommand\res{\mathop{\hbox{\vrule height 7pt width .3pt depth 0pt
\vrule height .3pt width 5pt depth 0pt}}\nolimits}

\title[Planar minimal graphs with analytic data]{Sharp boundary regularity for
planar minimal graphs\\ with analytic Dirichlet data}
\author{Camillo De Lellis}
\address{School of Mathematics, Institute for Advanced Study, Princeton, NJ 08540, USA\\
Gran Sasso Science Institute, L'Aquila, Italy}

\author{Guido De Philippis}
\address{Dipartimento di Matematica ``Tullio Levi-Civita'', Universit\`a di Padova, 35121 Padova, Italy}

\author{Alessio Figalli}
\address{Department of Mathematics, ETH Z\"urich, 8092 Z\"urich, Switzerland}
\date{}


\usepackage{aliascnt}

\newtheorem{thm}{Theorem}[section]

\newaliascnt{prop}{thm}
\newtheorem{prop}[prop]{Proposition}
\aliascntresetthe{prop}

\newaliascnt{lem}{thm}
\newtheorem{lem}[lem]{Lemma}
\aliascntresetthe{lem}

\newaliascnt{cor}{thm}
\newtheorem{cor}[cor]{Corollary}
\aliascntresetthe{cor}

\newaliascnt{claim}{thm}

\aliascntresetthe{claim}

\theoremstyle{definition}
\newaliascnt{defn}{thm}
\newtheorem{defn}[defn]{Definition}
\aliascntresetthe{defn}

\newaliascnt{notation}{thm}

\aliascntresetthe{notation}

\theoremstyle{remark}
\newaliascnt{remark}{thm}
\newtheorem{remark}[remark]{Remark}
\aliascntresetthe{remark}

\crefname{thm}{Theorem}{Theorems}
\crefname{prop}{Proposition}{Propositions}
\crefname{lem}{Lemma}{Lemmas}
\crefname{cor}{Corollary}{Corollaries}
\crefname{claim}{Claim}{Claims}
\crefname{defn}{Definition}{Definitions}
\crefname{notation}{Notation}{Notations}
\crefname{remark}{Remark}{Remarks}
\crefname{section}{Section}{Sections}
\crefname{subsection}{Section}{Sections}

\numberwithin{equation}{section}

\newcommand{\R}{\mathbb{R}}
\newcommand{\Ha}{\mathcal{H}}
\newcommand{\id}{\ensuremath{\mathds{1}}}
\newcommand{\spt}{\operatorname{spt}}
\newcommand{\Lip}{\operatorname{Lip}}
\newcommand{\dist}{\operatorname{dist}}

\begin{document}

\begin{abstract}
We study the boundary behavior of generalized (BV) solutions of the Dirichlet
problem for the minimal surface equation over a bounded set \(\Omega\). Such a generalized solution is always the graph of a smooth function \(v\) inside \(\Omega\), but might carry a vertical part over \(\partial\Omega\), which records the gap between the prescribed Dirichlet datum \(g\) and the boundary trace \(h\) attained by \(v\).
When \(\Omega\) is a bounded real-analytic planar domain and \(g\) is real analytic we
prove that \(h\) is Lipschitz, right and left differentiable at every point,
of class \(C^1\) on \(\partial\Omega\setminus\Sigma\), and
real analytic on \(\{h\neq g\}\setminus\Sigma\), where
\(\Sigma=\{\kappa=0\}\) is the zero set of the curvature of
\(\partial\Omega\). Also, at each point \(p\in\Sigma\cap\{h\neq g\}\), the two one-sided
branches extend real analytically across \(p\).

The proof recasts the vertical part near a contact point as a
scalar obstacle problem with a possibly degenerate right-hand side, whose model
case is \(\Delta v=H\,\id_{\{v>0\}}\) with \(H\ge0\) a polynomial.  These
results are optimal: for merely smooth data the trace may be discontinuous, while in the real-analytic case there are examples where corners do occur.
\end{abstract}

\maketitle

\section{Introduction}

Consider a bounded Lipschitz open set \(\Omega\subset \mathbb R^n\), a \(g\in L^1 (\partial \Omega)\), and let\footnote{The regularity of the functions allowed in the definition of \(\mathbf{m}\) does not play an important role: if \(\Omega\) and \(g\) are smooth, we can increase the regularity of \(u\) accordingly while \(\mathbf{m}\) stays the same.}
\begin{equation*}
\mathbf{m}:= \inf \left\{\int_\Omega \sqrt{1+|\nabla u|^2}: u|_{\partial \Omega} = g
\quad \mbox{and} \quad u\in W^{1,1} \cap C^1_{\mathrm{loc}} (\Omega)\right\}\, .
\end{equation*}
It is known, see e.g.~\cite[Sections 12.14 and 13.15]{Giusti}, that \(\mathbf{m}\) is not always attained, even when \(\partial \Omega\) and \(g\) are very smooth, because minimizing sequences can ``attach vertically'' to the cylinder \(\partial \Omega \times \mathbb R\): only some appropriate geometric conditions (most notably the nonnegativity of the mean curvature of \(\partial \Omega\)) guarantee attainment.

However, even when \(\mathbf{m}\) is not attained, we can still study the limiting behavior of any minimizing sequence \(\{v_k\}\), which is done as follows. Since the BV norm of \(v_k\) is bounded by the energy, we can extract an appropriate weak limit \(v\) that also belongs to \(BV\). Moreover, it can be shown (cf.~\cite[Theorem 14.13]{Giusti}) that \(v\) is then smooth, in fact real analytic, in the interior of \(\Omega\). From now on, we focus our attention on the boundary, and we start with the observation that, since \(v\) is a function of bounded variation, it has a well-defined trace on \(\partial \Omega\), which from now on is denoted by \(h\). By Gagliardo's Theorem we know that \(h\in L^1 (\partial \Omega)\). Moreover
\begin{equation}\label{e:attainment}
\mathbf{m} \geq \int_\Omega \sqrt{1+|\nabla v|^2}
\end{equation}
by lower semicontinuity of the area functional. However, as already mentioned, the inequality in~\eqref{e:attainment} might be strict because the graphs of the minimizing sequence \(\{v_k\}\) can cluster on some ``vertical portions'' on the cylinder \(\partial \Omega \times \mathbb R\).
To capture the latter behavior we:
\begin{itemize}
\item Fix some extension \(\tilde{g}\in C^1_{\mathrm{loc}}\cap W^{1,1} (\mathbb R^n\setminus \overline{\Omega})\) of \(g\) to the ``exterior'' of \(\Omega\);
\item Extend each \(v_k\) to \(\tilde{v}_k\) by setting it equal to \(\tilde{g}\) on \(\mathbb R^n \setminus \Omega\);
\item Consider the epigraphs \(E_k := \{(x,t)\in \mathbb R^{n+1} : t < \tilde{v}_k (x)\}\).
\end{itemize}
The sets \(E_k\) are sets of (locally) finite perimeter and their perimeter is uniformly bounded on any compact set of \(\mathbb R^{n+1}\). Up to subsequences we can extract a limit converging to a set of finite perimeter \(E\). The following theorem gives a good description of what happens when the inequality in~\eqref{e:attainment} is strict.

\begin{thm}[Giusti~{\cite[Theorem~14.6]{Giusti}}]\label{t:Giusti}
The set \(E\) defined above satisfies
\begin{align*}
E\cap (\Omega \times \mathbb R) &= \{(x,t): t<v(x)\}\\
E\setminus (\overline{\Omega} \times \mathbb R) &= \{(x,t) : t < \tilde{g} (x)\}\, .
\end{align*}
Moreover, if we denote by \(\partial^* E\) its reduced boundary (see e.g.~\cite[Definition~3.3]{Giusti}), then
\begin{equation*}
\mathbf{m} - \int_\Omega \sqrt{1+|\nabla v|^2} = \mathcal{H}^n (\partial^* E \cap (\partial \Omega\times \mathbb R)) = \int_{\partial \Omega} |g(x)-h(x)|\, .
\end{equation*}
\end{thm}

If the mean curvature of \(\partial \Omega\) is nonnegative in a neighborhood of \(x\in \partial \Omega\) and \(g\) is continuous at \(x\), then \(v\) and \(h\) are continuous at \(x\) and \(g(x)=h (x)\), cf.~\cite[Theorem~15.9]{Giusti}. Without further geometric restrictions, but under some mildly stronger assumptions on \(\Omega\) and \(g\), the following ``intrinsic regularity result'' was proved  by Simon and Bourni.

\begin{thm}[Simon~\cite{Simon1976}, Bourni~\cite{Bourni2011PMC}]\label{t:Bourni} Assume \(\partial \Omega\) and \(g\) are of class \(C^{1,\alpha}\) for some \(\alpha\in (0,1)\). Then
\(\partial^* E \cap (\overline{\Omega}\times \mathbb R)\) is an orientable \(C^{1,\alpha}\) submanifold of \(\mathbb R^{n+1}\) with boundary (in the classical sense of differential topology) and its boundary is the \((n-1)\)-dimensional oriented submanifold induced by the graph of \(g\) over \(\partial \Omega\).
\end{thm}

In~\cite{Simon1982}, Simon proved some further regularity of \(h\) but also showed examples of \(\Omega \subset \mathbb R^2\) of class \(C^k\) (with arbitrary \(k\in \mathbb N\)) and smooth boundary data \(g\) for which \(h\) has discontinuities. Later on Simon's example was extended by Williams to \(C^\infty\) domains in~\cite{Williams1986}.

In this paper we investigate what happens when \(\Omega\) is a two dimensional domain with real analytic boundary. Our main result shows that, when \(\Omega\subset \mathbb R^2\) is real analytic and \(g\) is real analytic, \(h\) is in fact Lipschitz, so that the ``vertical part'' of the generalized  minimizer is the closed region of \(\partial \Omega\times \mathbb R\) comprised between the two classical graphs of \(g\) and \(h\). Our analysis also yields some more refined regularity properties of \(h\). The main result is summarized in the following statement.

\begin{thm}\label{thm:main}
Let \(\Omega\subset\R^2\) be a bounded open set with real-analytic boundary,
let \(g\in C^\omega(\partial\Omega)\), and consider the Caccioppoli set \(E\), the function \(v\in BV (\Omega)\), and its trace \(h\) as in~\cref{t:Giusti}. Let
\(\kappa\) be the inward curvature of \(\partial\Omega\), and set \(\Sigma:=\{x\in\partial\Omega:\kappa(x)=0\}\).
Then:
\begin{enumerate}
\item[(i)] \(h\) is Lipschitz and admits right and left derivatives at every point. In particular the detachment set \(\{h\neq g\}\) is open and the support of the vertical part of \(\partial^*E\) over \(\partial\Omega\) is the closure of the region contained between the graphs of \(g\) and \(h\):
\[
\spt\bigl(\Ha^2\res (\partial^*E\cap(\partial\Omega\times\R))\bigr)
=\overline{\{(x,t):\min\{g(x),h(x)\}< t< \max\{g(x),h(x)\}\}}.
\]
\item[(ii)] \(h\in C^1(\partial\Omega\setminus\Sigma)\) and \(h\in C^\omega(\{h\neq g\}\setminus\Sigma)\).
\item[(iii)] If \(p\in\Sigma\cap\{h\neq g\}\), then its two one-sided branches
extend real analytically across \(p\).
\end{enumerate}
\end{thm}

\begin{remark}\label{remark-Sigma} In the above statement we are implicitly using the elementary fact that \(\Sigma\) is discrete (and thus finite): since \(\kappa\) is real analytic, if \(\Sigma\) were not discrete then it would contain a connected component \(\gamma\) of \(\partial \Omega\). Such a \(\gamma\) would be a straight line, contradicting the boundedness of \(\Omega\).
\end{remark}

\begin{remark}[Sharpness]
Both the conclusions and the assumptions are optimal: \cite{Williams1986} shows that even for \(C^\infty\) data
the trace \(h\) may be discontinuous, and we will discuss
in~\cref{rem:corner} that the two branches in~\cref{thm:main}(iii)
may meet at a corner even for real-analytic data.
 the relative boundary of \(\{h\ne g\}\), we expect the \(C^1\) regularity
of \(h\) outside \(\Sigma\) to be optimal in general, while at points in \(\Sigma\)
the two one-sided derivatives are distinct (see~\cref{r:Uraltseva}).
\end{remark}

\begin{remark}[Connection to the obstacle problem]
As in~\cite{LauLin85,Lin87}, in order to study regularity of the trace we reduce the Dirichlet problem to an obstacle problem in a suitable system of coordinates.
Fix a boundary point \(x\in \partial \Omega\) where \(g(x) \neq h(x)\) and to fix ideas assume that \(h<g\) in a neighborhood of \(x\). We change to new variables \(y\) so that \(T_x \partial \Omega \times \mathbb R\) is the horizontal plane of the new system of coordinates, while the vertical line \(\{(x,t): t\in \mathbb R\}\) is the axis \(\{(0, y_2, 0) : y_2 \in \mathbb R\}\). This change of variables recasts (locally) the ``vertical part'' as the contact set of a scalar obstacle problem for the gap \(w\) (in the \(y_3\) variable) between the surface and the obstacle \(\partial \Omega \times \mathbb R\).
This gap vanishes on the contact set and (since we are assuming that locally \(h<g\)) \(w\) solves a perturbation of
\[
  \Delta w = \bar{c} y_1^{2m}\,\id_{\{w>0\}},\qquad w\ge0,\qquad m \in \{0,1,2,3\ldots\}\, , \qquad \bar{c} > 0
\]
(the true equation is a
quasilinear minimal-surface analog).
The order \(2m\) and the positive coefficient reflect the fact that \(\kappa\) cannot change sign at \(x\), otherwise \(g(x)=h(x)\) (cf.~\cref{sec:sign-reduction}).  The fact that  \(v\)is a graph translates into the monotonicity of of \(w\) in the \(y_2\) direction: \(\partial_{y_2} w \ge0\), which  plays a crucial role in our arguments.
\end{remark}

The proof of~\cref{thm:main} is organized as follows.  In~\cref{sec:generalized-graphs} we recall the relaxed least-area formulation and
the boundary regularity input.~\Cref{sec:obstacle} reduces each one-sided
contact region to a scalar obstacle problem.~\Cref{sec:sign-reduction} shows
that the analytic forcing has even vanishing order and positive sign at every
genuine free-boundary point.  The interior and boundary free-boundary analyses
are then carried out in~\cref{sec:2d-blowups,sec:boundary-blowup} and assembled
in~\cref{sec:main-proof}, whose final part also
establishes the real analyticity of the two branches at the even-order corners.

\subsection*{AI declaration}
The main ideas and proofs of this paper were developed during the 2024–25 academic year while AF was visiting the Institute for Advanced Study in Princeton, as also noted in AF's Annual Director's Report submitted in June 2025.
 During the final stages of manuscript preparation, large language
models were used to assist with editing, checking some arguments, addressing minor technical details, and preparing the figures.

\subsection*{Acknowledgements}
The authors thank Emanuele Spadaro for drawing their attention to the question of trace regularity for minimal graphs.
GDP's research is  supported by the European Research Council (ERC) through Consolidator Grant No. 101169953 ``RISE''.\footnote{The views and opinions expressed are those of the authors only and do not necessarily reflect those of the European Union or the European Research Council.}
AF gratefully acknowledges the support of the Marvin V. and Beverly J. Mielke Fund, which made possible his stay at the Institute for Advanced Study (IAS) in Princeton during the 2024--25 academic year.

\section{Generalized graphs and boundary geometry}\label{sec:generalized-graphs}

From now on we assume that \(\Omega\) and \(g\) satisfy the assumptions of~\cref{thm:main}: several statements in what follows actually require much less regularity, but we leave the details to the reader. Observe that, since \(\Omega\subset\R^2\), its boundary \(\partial\Omega\) consists of finitely many embedded curves.  We will use \(E\), \(v\), \(h\) for the generalized
minimizing subgraph, its interior BV graph, and the boundary trace of the
introduction. By~\cref{t:Bourni}, the generalized graph
\[
G:=\partial^*E\cap(\overline\Omega\times\R)
\]
is an orientable \(C^{1,\alpha}\) \(2\)-dimensional submanifold with boundary given by the
graph of \(g\) over \(\partial\Omega\). Its vertical part \(\partial^* E \cap (\partial \Omega \times \mathbb R)\) is a closed subset of the ruled surface \(\partial \Omega \times \mathbb R\). The following elementary lemma is a simple consequence of~\cref{t:Bourni},~\cref{t:Giusti}, and the fact that \(E\cap (\Omega\times \mathbb R)\) is the epigraph of the locally smooth function \(v\).

\begin{lem}\label{l:elementary-contact}
For every \(x\in \partial \Omega\) the closed set \(G \cap (\{x\}\times \mathbb R)\) is either a singleton or a closed bounded vertical segment with \((x,g(x))\) as one endpoint. For \(\mathcal{H}^1\)-a.e. \(x\) the other endpoint is \((x,h(x))\).
\end{lem}

We can therefore define \(h (x)\) at \emph{every point} \(x\in \partial \Omega\) as follows:
\begin{itemize}
    \item \(h(x) = g(x)\) when \(G\cap (\{x\}\times \mathbb R)\) is a singleton;
    \item when \(G\cap (\{x\}\times \mathbb R)\) is not a singleton, \(h(x)\) is the only \(\tau\in \mathbb R\) such that \(G\cap (\{x\}\times \mathbb R)\) is the vertical segment joining \((x,g(x))\) and \((x,\tau)\).
\end{itemize}
We are now ready to decompose the vertical part of \(G\).

\begin{defn}\label{not:contact}
We define \(\Gamma := G \cap (\partial \Omega \times \mathbb R) = \partial^* E \cap (\partial \Omega \times \mathbb R)\)
and further decompose it into
\begin{align*}
\Gamma^+ &:=\{(x,t)\in\Gamma:t<g(x)\} = \{(x,t) : h(x)\leq t < g(x)\}\\
\Gamma^- &:=\{(x,t)\in\Gamma:t>g(x)\} = \{(x,t) : g(x) < t \leq h(x)\}\\
\Gamma_0 &:=\{(x,g(x)): x\in\partial\Omega,\ h(x)=g(x)\}.
\end{align*}
We furthermore subdivide \(\Gamma_0\) into
\begin{align*}
\Gamma_i &:= \Gamma_0 \cap \overline{\Gamma^+\cup \Gamma^-}\\
\Gamma_r &:= \Gamma_0 \setminus \Gamma_i =\{p\in\Gamma_0:B_\rho(p)\cap\Gamma\subset\Gamma_0
\text{ for some }\rho>0\}\, .
\end{align*}
\end{defn}

Recalling~\cite[Theorem~15.9]{Giusti}, the boundary datum is attained (i.e., \(h=g\)) in any neighborhood of a point \(x\in \partial\Omega\) where the inward
mean curvature is nonnegative. Since \(\Omega \subset \R^2\), the mean curvature of the boundary is just the classical curvature \(\kappa\). We summarize this result as follows:

\begin{prop}\label{prop:attachment}
If \(\kappa \geq 0\) in a neighborhood of \(x\in \partial \Omega\), then \(h=g\) in a neighborhood of \(x\) and thus \(p= (x,h(x)) \in \Gamma_r \subset \Gamma_0\).
\end{prop}

We next use~\cref{t:Bourni} to identify the tangent plane to the surface \(G\) at any point \(p= (x,t) \in \overline{\Gamma^+\cup\Gamma^-}\) and to prove that the two one-sided contact regions, \(\Gamma^+\) and \(\Gamma^-\), cannot meet, using the continuity of the conormal \(\nu_{\partial G} (p) \subset T_p G\) to the topological boundary \(\partial G\).

\begin{lem}\label{lem:vertical+separation}
If \(p=(x,t) \in \overline{\Gamma^+\cup\Gamma^-}\), then \(T_p G = T_p (\partial \Omega \times \mathbb R) = (T_x \partial \Omega) \oplus \mathbb R e_3\). Moreover \(\overline{\Gamma^+}\cap\overline{\Gamma^-}=\emptyset\). In particular, for every \(p\in\overline{\Gamma^+\cup\Gamma^-}\) there is \(\rho>0\) such that \(B_\rho(p)\) meets only one of \(\Gamma^+\) and \(\Gamma^-\).
\end{lem}
\begin{proof}
In order to prove the first statement observe first that, according to~\cref{t:Bourni}, when \(p=(x,t)\in \Gamma^+\cup \Gamma^-\) the point \(p\) is \emph{not} contained in the differential topological boundary \(\partial G\) of the surface \(G\). In particular \(B_\rho (p) \cap G\) must be a topological disk containing \(p\) in its interior for every sufficiently small \(\rho\). On the other hand if \(T_p G\) did not coincide with \(T_p (\partial \Omega \times \mathbb R)\), it would need to be transversal to it, and given that \(G\) is entirely contained in the closed set \(\overline{\Omega} \times \mathbb R\), \(B_\rho (p)\cap G\) would be a topological half-disk containing \(p\) in its boundary, for every sufficiently small \(\rho\). For the remaining points \(p\in\overline{\Gamma^+\cup\Gamma^-}\), the claim follows by approximation and continuity of the tangent plane.

\medskip

The last statement is an obvious consequence of the second and we will therefore focus on the latter.
We let \(\nu_{\partial G}\) be the exterior unit conormal to the oriented \(C^{1,\alpha}\) surface \(\partial G\), where \(\partial G\) is understood as the boundary, in the sense of differential topology, of \(G\). Consider \(p= (x,t)\in \Gamma^+\) and let \(\bar{p} = (x, g(x))\). From the first statement we know that \(\nu_{\partial G} (\bar{p})\in T_x \partial \Omega \oplus \mathbb R e_3\). Since \(h(x)< g(x)\), the points \(\{(x,s): h(x)\leq s< g(x)\}\) belong to \(G\): given that the conormal must be pointing ``outward'' we must have \(\nu_{\partial G} (\bar{p})\cdot e_3 \geq 0\). Furthermore, since \(\nu_{\partial G} (\bar{p})\) is orthogonal to the graph of \(g\), the latter inequality upgrades to
\(\nu_{\partial G} (\bar{p}) \cdot e_3 \geq (1+\Lip(g)^2)^{-1/2}\). However when \(p= (x,t)\in \Gamma^-\) the very same considerations give \(\nu_{\partial G} (\bar{p}) \cdot e_3 \leq - (1+\Lip(g)^2)^{-1/2}\). But then the continuity of \(\nu_{\partial G}\) prevents the existence of a point where both \(\Gamma^+\) and \(\Gamma^-\) accumulate.
\end{proof}

\section{Reduction to an obstacle problem}\label{sec:obstacle}

Fix \(p=(x,g(x))\in\overline{\Gamma^+\cup\Gamma^-}\), that is, either \(p\in\Gamma_i\) or \(h(x)\neq g(x)\).  By~\cref{lem:vertical+separation}, after shrinking the
neighborhood of \(p\), only one of \(\Gamma^+\) and \(\Gamma^-\) accumulates at \(p\). To fix ideas assume it is a point of accumulation of \(\Gamma^+\): the other case is treated analogously, after reversing the
vertical orientation.

\begin{figure}[hbtp]
\begin{tikzpicture}
    \begin{axis}[
        view={45}{35},
        axis lines=center,
        width=12cm,
        height=10cm,
        plot box ratio=1 2.5 1,
        xmin=-0.6, xmax=0.6,
        ymin=-1.5, ymax=1.5,
        zmin=-0.1, zmax=0.3,    %
        xlabel={\(y_1\)}, ylabel={\(y_2\)}, zlabel={\(y_3\)},
        grid=major,
        xticklabels=\empty,
        yticklabels=\empty,
        zticklabels=\empty,
        colormap={lightmap}{color(0cm)=(white); color(1cm)=(yellow!20); color(2cm)=(orange!20)}
    ]

    \def\R{3.0}
    \def\th{11}
    \def\yheight{1.5}

    \addplot3 [
        surf,
        colormap name=lightmap,
        shader=faceted,
        opacity=0.6,
        faceted color=orange!40,
        domain=-\th:\th,
        domain y=-\yheight:\yheight,
        samples=25, samples y=10
    ] ({\R*sin(x)}, {y}, {\R*(cos(x)-1)});

    \addplot3 [
        surf,
        color=blue!10,
        opacity=0.3,
        dashed,
        domain=-0.6:0.6,
        domain y=-1.5:1.5,
        samples=2
    ] {0};

    \draw[red!60, ultra thick] (axis cs: 0, -\yheight, 0) -- (axis cs: 0, \yheight, 0);

    \end{axis}
\end{tikzpicture}

\caption{The set \(\partial \Omega \times \mathbb R\) as a (local) graph of a function \(\varphi (y_1, y_2)\) in new coordinates \(y_1, y_2, y_3\). The \(y_2\) axis coincides with the \(x_3\) axis and thus \(\varphi\) depends only on \(y_1\). The \(y_3\) axis is orthogonal and points inward to \(\Omega\). Therefore in the picture we are looking at a portion where the curvature \(\kappa\) of \(\partial \Omega\) is negative, which are the points of interest in our analysis.}\label{f:phi}
\end{figure}
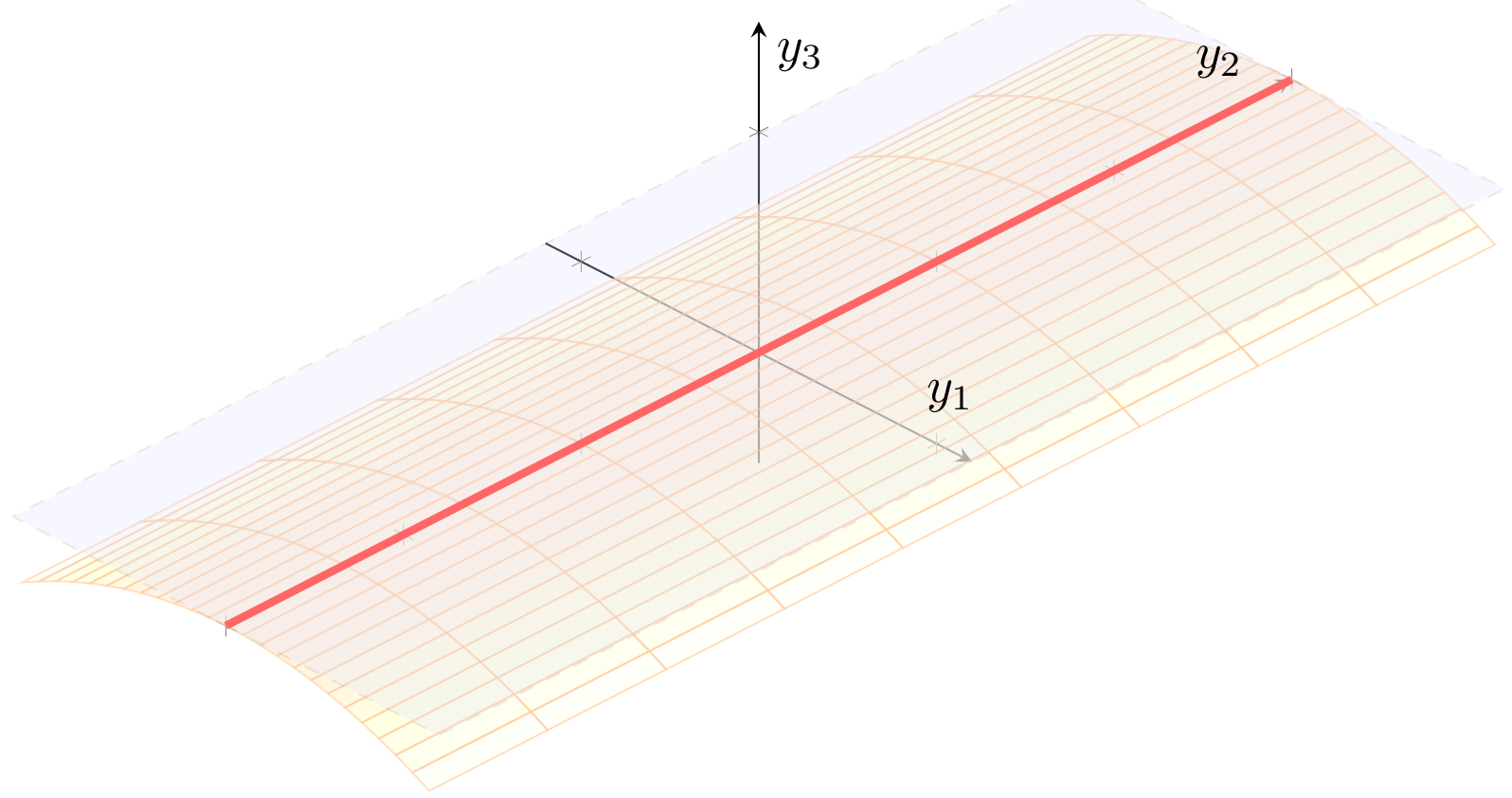

We start by selecting a new system of coordinates
\((y_1,y_2,y_3)\) so that:
\begin{itemize}
\item \(p = (0, a, 0)\) for some \(a\in \mathbb R\),
\item \(T_x\partial\Omega\) is the axis \(\{(y_1,a,0): y_1\in \mathbb R\}\),
\item \(\{x\}\times\R\) is the axis \(\{(0,y_2,0): y_2\in \mathbb R\}\),
\item and \((0,0,1)\) is the inward normal to
the cylinder \(\Omega \times \mathbb R\).
\end{itemize}
Furthermore we select a sufficiently small neighborhood \(\mathcal{U} = B_\rho (x) \times \mathbb R \subset \mathbb R^2\times \mathbb R\) of the line \(\{x\}\times \mathbb R\) so that \(\partial \Omega \cap \mathcal{U}\) is graphical over the plane \(\pi_o := T_x \partial \Omega \oplus \mathbb R e_3 = \{y : y_3=0\}\). \(\mathcal{U} \cap (\partial\Omega\times\R)\) is in turn described as a graph \(\{(y_1, y_2, \varphi (y_1)): |y_1|< \rho_0, y_2\in \mathbb R\}\), where \(\varphi\) is a real-analytic function over some interval \((-\rho_0, \rho_0)\), cf.~\cref{f:phi}.

We consider the open subset \(\Xi\) of \(\pi_o\) given by the orthogonal projection onto \(\pi_o\) of the epigraph of \(g\) in \(\partial \Omega \cap \mathcal{U}\). The ``bottom part'' \(\partial_b \Xi\) of the boundary of \(\Xi\) is the projection onto \(\pi_o\) of the portion of the graph of \(g\) (a \(1\)-dimensional submanifold of the cylinder \(\partial \Omega \times \mathbb R\)) which lies in \(\mathcal{U}\). We can describe it as
\[
\partial_b \Xi = \{y_2=\tilde g(y_1) : |y_1| < \rho_0\}
\]
for some \(\tilde g\in C^\omega (-\rho_0, \rho_0)\).
We now choose the sign of the coordinate \(y_2\) so that
\[
\Xi=\{(y_1,y_2):y_2>\tilde g(y_1)\}.
\]
Finally, thanks to~\cref{t:Bourni} and~\cref{lem:vertical+separation}, in a neighborhood of \(p\)
the surface \(G\) can be described as the graph of a \(C^{1,\alpha}\) function \(\tilde u\) over \(\Xi \cap V\), for some neighborhood \(V\) of the projection of \(p\) onto \(\pi_o\), lying above the function \(\varphi\). In fact such a description is valid also in a neighborhood of any point \((y_1,y_2,\varphi (y_1))\) in \(\mathcal{U}\cap \Gamma^+\).

\begin{lem}\label{lem:obstacle}
Upon choosing \(\rho\) (and hence \(\rho_0\)) sufficiently small, near any \(q\in \partial_b \Xi\) and near any \(q = (y_1, y_2, \varphi (y_1)) \in \Gamma^+\cap \mathcal{U}\), the generalized surface \(G\) is the graph of a function \(\tilde u\), on an open subset of \(\Xi\), with
\begin{enumerate}
\item \(\tilde u\ge\varphi\);
\item \(\tilde u=\varphi\) on \(\partial_b \Xi\);
\item \(\partial_{y_2}\tilde u\ge0 \).
\end{enumerate}
\end{lem}

\begin{proof}
 Properties (1) and (2) are immediate from the construction.
For (3), \(E\) is a subgraph in the ``original'' vertical direction \(t\): if \((x,t)\in E\) and \(s\ge0\) then \((x,t-s)\in E\).  In the \(y\)-variable we have that, near \(p\), \(E=\{y_3<\tilde u(y_1,y_2)\}\), and the subgraph property becomes
\(y_3<\tilde u(y_1,y_2)\), implying therefore that \(y_3<\tilde u(y_1,y_2+s)\) for \(s\ge0\).
This proves that \(\partial_{y_2}\tilde u\ge0\), as desired.
\end{proof}

We now choose a set \(\Lambda\subset \pi_o\) which includes all the points \(q\) of~\cref{lem:obstacle} and over which a function \(\tilde{u}\) as in~\cref{lem:obstacle} can be globally defined. It is convenient to assume that \(\Lambda\) has the following property:
\begin{itemize}
    \item[(I)] for all \(s\in (-\rho_0, \rho_0)\) the set \(\{y_1 =s\}\cap \Lambda\) is a half-open interval which contains its lowest endpoint \((s, \tilde{g} (s))\).
\end{itemize}
This can be easily achieved as follows. For every \(s\in (-\rho_0, \rho_0)\) we let \(\tau\geq \tilde{g} (s)\) be the maximal number such that \((s, \tau, \varphi (s))\) belongs to \(\Gamma^+\cup \Gamma_0\) and then choose a set \(V_s\) of the form
\[
V_s = \{(y_1, y_2) : |y_1-s|< \eta, \quad \tilde{g} (y_1) \leq y_2 < \tau +\eta\}
\]
over which a function \(\tilde{u}\) as above can be defined. This is possible for a sufficiently small \(\eta\), where \(\eta\) might depend on \(s\). We then take \(\Lambda\) to be the union of all \(V_s\) as \(s\) varies in \((-\rho_0, \rho_0)\): it is easy to see that \(\Lambda\) satisfies (I). Moreover, from now on \(\partial_b \Xi = \{(y_1, \tilde{g} (y_1)) : -\rho_0 < y_1 < \rho_0\}\) will be denoted by \(\partial_b \Lambda\), as it is indeed the bottom part of the boundary of \(\Lambda\).

\begin{figure}[hbtp]
\begin{tikzpicture}

  \begin{axis}[
    axis lines=middle,
    xlabel=\(y_1\),
    ylabel=\(y_2\),
    xtick={-4, 4},
    xticklabels={{\((-\rho_0, 0)\)}, {\((\rho_0, 0)\)}},
    ytick=\empty,
    label style={font=\small},
    tick label style={font=\footnotesize},
    every axis label/.append style={black},
    axis line style={black},
    xmin=-5.5, xmax=5.5,
    ymin=-1.0, ymax=8,
    samples=200,
    unit vector ratio*=1 1,
    clip=false
  ]

    \pgfmathsetmacro{\RVAL}{4}
    \pgfmathsetmacro{\acon}{2}

    \addplot [name path=G, thick, black, domain=-\RVAL:\RVAL] {0.05*x^2 + 1};

    \addplot [name path=H, thick, black, domain=-\RVAL:\RVAL] {0.05*x^2 + 1 + 0.3*max(0, -\acon-x) + 1.2*max(0, x-\acon)};

    \addplot [black!15] fill between[of=G and H];

    \pgfmathsetmacro{\gval}{0.05*\RVAL*\RVAL + 1}

    \draw [black, thin] (axis cs: \RVAL, \gval) -- (axis cs: \RVAL, 7.5);
    \draw [black, thin] (axis cs: -\RVAL, \gval) -- (axis cs: -\RVAL, 7.5);

    \node[font=\large] at (axis cs: 1.5, 3.5) {\(\Xi\)};

    \node[font=\footnotesize] at (axis cs: 1.5, 0.75) {\(\partial_b \Xi\)};

  \end{axis}

\end{tikzpicture}
\caption{The graphs of the functions \(\tilde{g}\) and \(\tilde{h}\): the shaded region is the projection onto \(\pi_o\) of the set \(\overline{\Gamma^+}\). The function \(\tilde{u}\) is defined over a neighborhood of this set and describes, locally, \(G = \partial^* E\): the shaded region is therefore the contact set \(\{\tilde{u} =\varphi\}\).}\label{f:contact-set}
\end{figure}
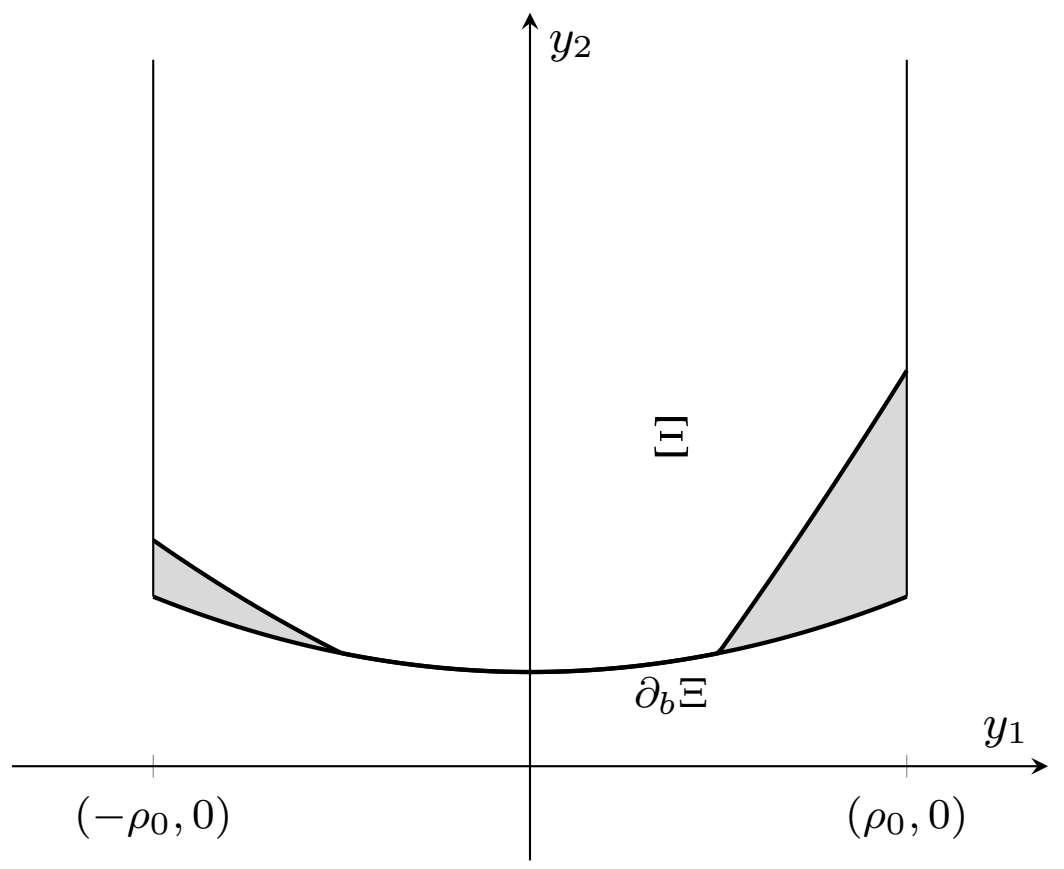

In the \(y\)-coordinates,\(\overline{\Gamma^+}\cup \Gamma_0\) projects onto the coincidence set
\(\{\tilde u=\varphi\}\cap\Lambda\).  By~\cref{lem:obstacle},  this set is given by
\(\{\tilde g(y_1)\le y_2\le\tilde h(y_1)\}\) for some function
\(\tilde h(y_1)\ge\tilde g(y_1)\), so the free boundary
\(\partial\{\tilde u>\varphi\}\cap\Lambda\) is the graph
\(\{y_2=\tilde h(y_1)\}\).  Just as \(\tilde g\) corresponds to \(g\) in these
coordinates, \(\tilde h\) corresponds to the projection onto \(\pi_o\) of the graph of the trace \(h\); in particular, the regularity of \(h\) is that of \(\tilde h\); cf.~\cref{f:contact-set}.

Next, from the area minimality of \(G\) we derive that
\(\tilde u\)  solves an obstacle problem with obstacle
\(\varphi\). The idea is to use competitors for \(\tilde u\) which are functions \(w\) with
\(w\ge\varphi\) and \(w-\tilde u\in W^{1,1}_0 (\Lambda)\): replacing the corresponding piece of
\(G\) by the graph of such a \(w\) gives the variational inequality
\begin{equation}\label{e:variational-ineq}
\int D_pF(\nabla \tilde u)\cdot \nabla(w-\tilde u)\ge0,\qquad F(p):=\sqrt{1+|p|^2}\, \, .
\end{equation}
We summarize the statement in the following lemma and we provide the details of the proof for the reader's convenience.

\begin{lem}\label{lem:variational-ineq}
The variational inequality~\eqref{e:variational-ineq} holds for every \(w\geq \varphi\) with \(w-\tilde{u}\in W^{1,1}_0 (\Lambda)\).
\end{lem}

\begin{proof}
First of all consider a function \(w\geq \varphi\) such that \(w-\tilde{u}\in C^1_c (\Lambda)\). For every \(\varepsilon\geq 0\) sufficiently small we let \(w_\varepsilon = \tilde{u} + \varepsilon (w-\tilde{u})\). If we choose a sufficiently small neighborhood \(V\subset \pi_o\) of the support of \(w-\tilde{u}\) and a suitable \(\delta > \|\tilde{u}\|_{C^0}\), in the open set \(W:= \{(y_1, y_2, y_3): (y_1, y_2)\in V, |y_3|< \delta\}\) the set \(E\) is given by the epigraph of \(\tilde{u}\):
\[
E\cap W = \{ (y_1, y_2, y_3): (y_1, y_2)\in V, -\delta < y_3 < \tilde{u} (y_1, y_2)\}\, ,
\]
(recall that we are dealing with the \(\Gamma^+\) case: for the \(\Gamma^-\) case the considerations are entirely analogous).
We now consider the set \(E'\) given by the union of \(E\setminus W\) and the epigraph of \(w_\varepsilon\) inside \(W\), provided \(\varepsilon\) is small enough to ensure \(\|w_\varepsilon\|_{C^0} < \delta\). Since \(w_\varepsilon \geq \varphi\), this set \(E'\) coincides with \(E\) outside the cylinder \(\Omega \times \mathbb R\). In particular we can apply~\cite[Lemma~14.7]{Giusti} to conclude that
\[
f(\varepsilon) := \mathcal{H}^n (\partial^* E' \cap (\overline{\Omega}\times \mathbb R)) - \mathcal{H}^n (\partial^* E \cap (\overline{\Omega}\times \mathbb R)) \geq 0
\, .
\]
In turn, using the \(C^1\) regularity of \(w_\varepsilon\) and \(\tilde{u}\) we compute
\[
f (\varepsilon) = \int \left(\sqrt{1+ |\nabla w_\varepsilon|^2} - \sqrt{1+|\nabla \tilde{u}|^2}\right)\, .
\]
Differentiating in \(\varepsilon\) at \(0^+\) we get \(f'(0) \geq 0\), which is precisely the inequality~\eqref{e:variational-ineq}. A standard approximation procedure allows us to conclude the same inequality in the general case of \(w-\tilde{u} \in W^{1,1}_0 (\Lambda)\).
\end{proof}

Next, for a \(C^{1,1}\) graph \(f\), we compute its mean curvature as
\[
H_f=(1+|\nabla f|^2)^{-1/2}
\left(I-\frac{\nabla f\otimes \nabla f}{1+|\nabla f|^2}\right):D^2f
\]
and introduce the symmetric matrix
\[
A(\nabla f) := (1+|\nabla f|^2)^{-1/2}
\left(I-\frac{\nabla f\otimes \nabla f}{1+|\nabla f|^2}\right)\, ,
\]
so that \(H_f = A (\nabla f) : D^2 f\).

Now, set \(\bar u:=\tilde u-\varphi\).
The next proposition proves that \(\bar u\in C^{1,1}_{\mathrm{loc}}\) and solves the obstacle-type problem
\begin{equation}\label{eq:ubar-pde}
A(\nabla\varphi+\nabla\bar u):D^2\bar u+
D^2\varphi:\bigl[A(\nabla\varphi+\nabla\bar u)-A(\nabla\varphi)\bigr]
=H\,\id_{\{\bar u>0\}},\qquad H:=-H_\varphi .
\end{equation}
Since \(\varphi\) depends only on \(y_1\), we have \(H_\varphi=\varphi''/(1+\varphi'^2)^{3/2}\), which is the curvature of \(\partial\Omega\) with respect to the inward normal \(e_3\) at the point with coordinate \(y_1\). Hence \(H=-\kappa\) in these coordinates: \(H>0\) exactly where \(\partial\Omega\) is strictly concave, and \(H<0\) where \(\kappa>0\), which is where~\cref{prop:attachment} applies.

\begin{prop}\label{prop:c11}
The function \(\bar u\) is locally \(C^{1,1}\) and solves~\eqref{eq:ubar-pde}.
\end{prop}

\begin{proof}
This proof is classical and we report it for the reader's convenience.
Inside \(\{\tilde u>\varphi\}\) both signs of variations are admissible in the
variational inequality, so the Euler--Lagrange equation gives \(H_{\tilde u}=0\);
in addition, upward variations from the obstacle give \(H_\varphi\le0\) on the contact
set. Classical theory for variational inequalities tells us that \(\tilde u\in W^{2,2}_{\mathrm{loc}}\) (cf.~\cite[Theorems~2.1 and~2.2]{KinderlehrerStampacchia1980}). Thus, since \(\nabla\tilde u=\nabla\varphi\) on \(\{\tilde u=\varphi\}\), we conclude that
\(D^2\tilde u=D^2\varphi\) a.e. on the very same set. In particular we get
\[H_{\tilde u}-H_\varphi=-H_\varphi\,\id_{\{\tilde u>\varphi\}}.\]
Note now that
\[
H_{\tilde u}-H_\varphi
=A(\nabla\varphi+\nabla\bar u):D^2\bar u
+[A(\nabla\varphi+\nabla\bar u)-A(\nabla\varphi)]:D^2\varphi.
\]
Also, since \(\varphi\) is smooth and \(\bar u\in C^{1,\alpha}\) (by~\cref{t:Bourni}),
it follows that
the matrix \(A(\nabla\varphi+\nabla\bar u)\) is uniformly elliptic and H\"older continuous.
Thus, the optimal \(C^{1,1}\) regularity  follows
from the classical obstacle problem
(see~\cite{Caffarelli1977,Caffarelli1998,PetrosyanShahgholianUraltseva2012}) and, for instance, the
quasilinear area-type extensions in~\cite{FocardiGeraciSpadaro2017}.
\end{proof}

By~\cref{lem:obstacle}, \(\partial_{y_2}\bar u\ge0\); the next proposition upgrades
this to strict positivity on \(\{\bar u>0\}\).

\begin{prop}\label{prop:strict-monotonicity}
On \(\{\bar u>0\}\) we have \(\partial_{y_2}\bar u>0\).
\end{prop}

\begin{proof}
We begin by noticing that, inside \(\{\bar u >0\}\), the function \(\tilde u\) solves the equation \(H_{\tilde u}=0\). Also, since \(\varphi\) is independent
of \(y_2\), \(\chi:=\partial_{y_2}\bar u=\partial_{y_2}\tilde u\).  Differentiating the PDE
\(A_{ij}(\nabla\tilde u)\partial_{ij}\tilde u=0\) gives
\begin{equation}\label{e:PDE-for-u_2}
A_{ij}(\nabla\tilde u)\partial_{ij}\chi
+\partial_{p_k}A_{ij}(\nabla\tilde u)\,\partial_{ij}\tilde u\,\partial_k \chi=0 .
\end{equation}
Since \(\chi \geq 0\) (by~\cref{lem:obstacle}(3)), the strong
maximum principle implies either \(\chi>0\) or \(\chi\equiv0\) on any connected component \(\mathcal{V}\) of \(\{\bar u >0\}\). But \(\chi \equiv 0\) would mean that \(\partial_{y_2} \bar u\) vanishes everywhere on the connected component, which is impossible. Indeed recall that for every \(s\in (-\rho_0, \rho_0)\) the intersection \(\{y_1=s\}\cap \Lambda\) is an interval \(\{(s,t) : \tilde{g} (s) < t < \tau\}\), while the subset of it where \(\bar{u}\) vanishes is either a closed interval with bottom endpoint \((s, \tilde{g} (s))\), or the singleton \((s, \tilde{g} (s))\). We conclude that the remaining positivity interval is either all contained in \(\mathcal{V}\) or it does not intersect \(\mathcal{V}\). But if it is contained in \(\mathcal{V}\) then \(\bar u\) vanishes identically on it and, by the arbitrariness of \(s\), \(\bar u\) must vanish identically on \(\mathcal{V}\).
\end{proof}

\section{Even order and positive sign of the forcing}\label{sec:sign-reduction}

In this section we will start exploiting the fact that our boundary data are analytic. We fix our attention on a point \(q\in \Gamma\) which is an ``interior free-boundary point'', namely
\(q= (x, h(x))\) with \(x\in \partial \Omega\) and \(h(x) \neq g (x)\). We next use the coordinates just introduced in~\cref{sec:obstacle}, where the point \(p\in\overline{\Gamma^+\cup\Gamma^-}\) is chosen to be \((x,g(x))\). In particular \(\tilde{h} (0) > \tilde{g} (0)\).
Recall also that \(\varphi(0)=0\) and \(\nabla\varphi(0)=0\). Hence, given that \(\varphi\) is analytic,
\[
H (y_1)=c\,y_1^{k}+O(y_1^{k+1}),\qquad c\ne 0\, ,
\]
for some integer \(k\ge0\) (recall that \(H\) cannot vanish identically because of~\cref{remark-Sigma}).
The purpose of this section is to show that the function \(H\) takes the form
\[
H (y_1)=c\,y_1^{2m}+O(y_1^{2m+1}),\qquad c>0.
\]
We begin with a Hopf-lemma continuity statement: if \(\tilde h=\tilde g\) on one
side of a boundary point, then \(\tilde h=\tilde g\) at the
point itself.

\begin{lem}\label{lem:boundary-empty-side}
 Assume that \(\tilde h=\tilde g\) on at least one of the intervals
\((-\rho,0)\) and \((0,\rho)\), for some \(\rho>0\). Then \(\tilde h(0)=\tilde g(0)\).
\end{lem}

\begin{proof}
Without loss of generality, we can assume that \(\tilde h=\tilde g\)
on \((-\rho,0)\).
Assume by contradiction \(\tilde h(0)>\tilde g(0)\), consider the segment \(I_0:=\{0\}\times (\tilde g(0),\tilde h(0))\), and define
\[
s_0:=\frac{\tilde h(0)+\tilde g(0)}2,\qquad r_0:=\frac{\tilde h(0)-\tilde g(0)}2 .
\]
Since \(\tilde h=\tilde g\)
on \((-\rho,0)\), we can find \(r \in (0,r_0)\) small such that the two-dimensional ball \(B_r(-r,s_0)\) is contained inside \(\mathcal{O}=\{\bar u>0\}\) and touches \(I_0\subset \partial \mathcal{O}\) at the point \((0,s_0)\), cf.~\cref{f:Hopf-1}.

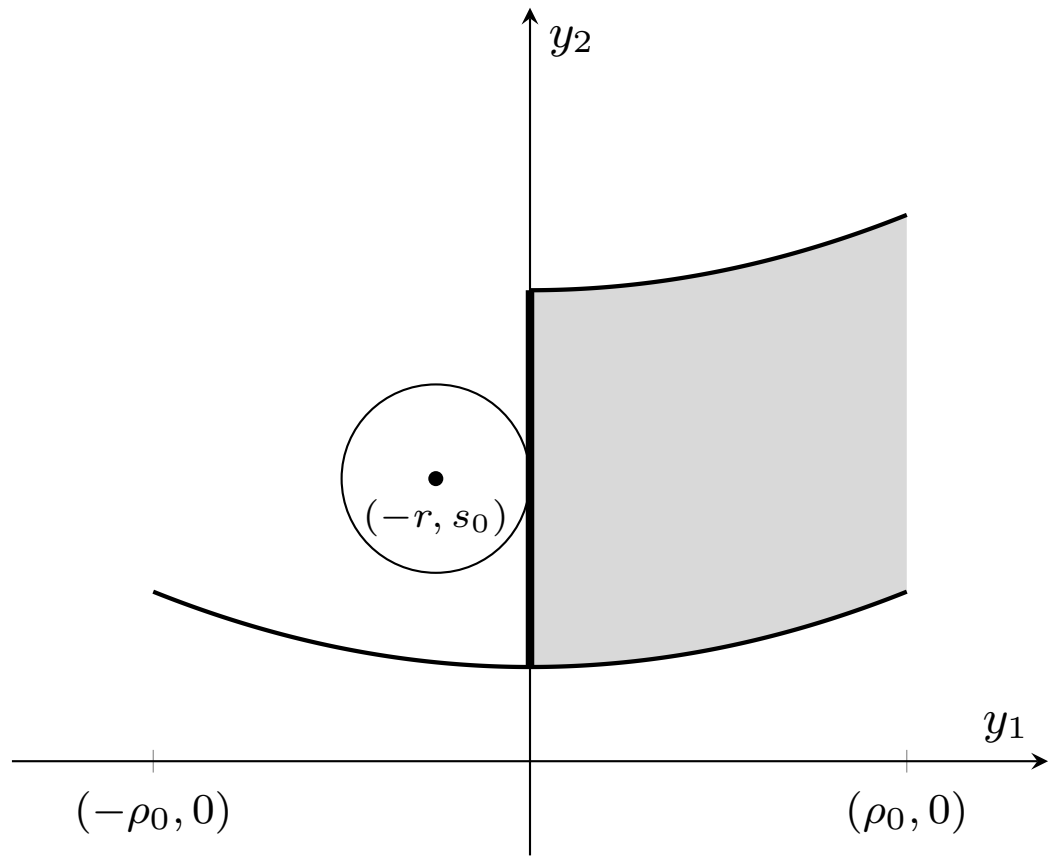
\begin{figure}
\begin{tikzpicture}

  \begin{axis}[
    axis lines=middle,
    xlabel=\(y_1\),
    ylabel=\(y_2\),
    xtick={-4, 4},
    xticklabels={{\((-\rho_0, 0)\)}, {\((\rho_0, 0)\)}},
    ytick=\empty,
    label style={font=\small},
    tick label style={font=\footnotesize},
    every axis label/.append style={black},
    axis line style={black},
    xmin=-5.5, xmax=5.5,
    ymin=-1.0, ymax=8,
    samples=200,
    unit vector ratio*=1 1,
    clip=false
  ]

    \pgfmathsetmacro{\RVAL}{4}

    \addplot [name path=G, thick, black, domain=-\RVAL:\RVAL] {0.05*x^2 + 1};

    \addplot [name path=H, thick, black, domain=0:\RVAL] {0.05*x^2 + 5};

    \addplot [black!15] fill between[of=H and G, soft clip={domain=0:\RVAL}];

    \draw [black, ultra thick] (0, 1) -- (0, 5);

    \draw [black, thin] (-1, 3) circle [radius=1];

    \node [circle, fill=black, inner sep=1pt] at (-1, 3) {};
    \node [below=2pt, font=\scriptsize] at (-1, 3.1) {\((-r, s_0)\)};

  \end{axis}
\end{tikzpicture}
\caption{The disk \(B_r (-r, s_0)\). We apply Hopf's lemma at the point \((0, s_0)\) where the boundary touches the vertical segment \(I_0\subset \{\bar{u} = 0\}\). The disk is chosen so that its interior is contained in \(\{\bar{u}>0\}\).}\label{f:Hopf-1}
\end{figure}

Set \(\chi:=\partial_{y_2}\bar u\).  By~\cref{prop:strict-monotonicity}, \(\chi>0\) in
\(\mathcal{O}\); also, since \(\nabla \bar u\) vanishes on the free boundary, \(\chi=0\) on \(I_0\). Recall that \(\chi\) satisfies~\eqref{e:PDE-for-u_2}, which we can rewrite as the linear PDE
\[
a_{ij}(y)\partial_{ij}\chi+b_i(y)\partial_i \chi=0 \qquad \text{inside \(\mathcal{O}\),}
\]
with \(a_{ij}\) uniformly elliptic. So Hopf's lemma yields
\begin{equation}\label{e:Hopf-100}
-\partial_{y_1} \chi(0,s_0)>0\, ,
\end{equation}
where the left-hand side is interpreted as the normal derivative of the left trace of \(\chi\) over \(\partial B_r (-r, s_0)\).
On the other hand, recall that
\(\nabla \bar u\) vanishes on the free boundary, hence it vanishes on \(I_0\). We would like to differentiate this equality along \(I_0\) to deduce that \(\partial_{y_2} \partial_{y_1} \bar u=0\) on \(I_0\) and find a contradiction with~\eqref{e:Hopf-100}. This requires some care, because \(\bar u\) is only \(C^{1,1}\). Note however that the partial derivative \(\partial_{y_1} \bar u\) solves the linearized minimal-surface equation in the half-ball \(B^-:= \{y_1<0\}\cap B_{r_0}(0,s_0)\). Hence, since \(\partial_{y_1} \bar u\) vanishes on \(I_0=\{y_1=0\}\cap B_{r_0}(0,s_0)\) and the coefficients of the PDE are H\"older continuous, elliptic regularity tells us that \(\partial_{y_1} \bar u\) extends to a \(C^{1,\alpha}\) function on the closure of \(B^-\). In particular we can compute the tangential derivative of \(\partial_{y_1} \bar u\) along the interval \(I_0\subset \partial B^-\) at \((0,s_0)\) and we get
\[
0=\partial_{y_2}\partial_{y_1}\bar u(0,s_0)=\partial_{y_1}\partial_{y_2}\bar u(0,s_0)=\partial_{y_1}\chi(0,s_0)<0,
\]
a contradiction. This proves that \(\tilde h(0)=\tilde g(0)\).
\end{proof}

Using the previous lemma, we prove that at an interior free-boundary point, the first nonzero analytic
term of \(H\) must have even order and positive coefficient.

\begin{lem}\label{lem:no-sign-change}
Assume that \(\tilde h(0)>\tilde g(0)\). Then \(H (y_1)=\bar c y_1^{2m}+O(y_1^{2m+1})\) for some \(m\geq 0\) and \(\bar c>0\).
\end{lem}

\begin{proof}
By analyticity, \(H(y_1)=\bar c y_1^{k}+O(y_1^{k+1})\) for some \(k \geq 0\) and \(\bar c \neq 0\).
If, by contradiction,  either \(k\) is odd or \(\bar c<0\), we can find \(\rho>0\) such that \(H<0\) on at least one of the intervals
\((-\rho,0)\) and \((0,\rho)\). Without loss of generality, we can assume \(H<0\) on \((-\rho,0)\). Then~\cref{prop:attachment} implies that
\(\tilde h=\tilde g\) on \((-\rho,0)\), and therefore~\cref{lem:boundary-empty-side}
gives \(\tilde h(0)=\tilde g(0)\), a contradiction.
\end{proof}

\section{Interior blow-ups}\label{sec:2d-blowups}

In this section we complete the study of ``interior'' free-boundary points. We assume the same setup as in the previous section, in particular the free boundary point is assumed to be the origin and moreover \(\tilde h(0) > \tilde g(0)\), while \(\varphi(0)=0\) and \(\nabla\varphi (0)=0\).
In addition, thanks to~\cref{lem:no-sign-change},
\begin{equation*}
H(y_1)=\bar c\,y_1^{2m}+O(y_1^{2m+1}) \qquad \mbox{for some \(\bar c>0\)}.
\end{equation*}
Our first goal is to prove the monotonicity formulas that detect the homogeneity
\(\lambda=2m+2\) of the blow-ups.
To this end, we first prove the following:

\begin{lem}\label{lem:finite-order-upper-detachment}
Set \(\lambda=2m+2\).  There are \(C,r_0>0\) such that, for \(0<r<r_0\),
\[
\sup_{B_r}\bar u\le Cr^\lambda,\qquad
\sup_{B_r}|\nabla\bar u|\le Cr^{\lambda-1},\qquad
\sup_{B_{r/2}}|D^2\bar u|\le Cr^{\lambda-2}.
\]
\end{lem}

\begin{proof}
Since \(0\) is a free-boundary point,
\[
\bar u(0)=0,\qquad \nabla\bar u(0)=0 .
\]
The free-boundary \(C^{1,1}\) estimate for obstacle-type variational
inequalities gives, with a constant independent of \(\|\bar u\|_{L^\infty}\),
\[
\sup_{B_r}|D^2\bar u|
\le C\|H\|_{L^\infty(B_{2r})}.
\]
Indeed, this follows from the usual growth estimate at a free-boundary
point, followed by interior \(C^{1,1}\) estimates for the elliptic PDE (in our case the minimal surface equation) in the positivity set, cf.~\cite{FocardiGeraciSpadaro2017}.  Since
\(\|H\|_{L^\infty(B_{2r})}\le Cr^{\lambda-2}\), the Hessian estimate follows.
Taylor's formula at the free-boundary point then gives the estimates for
\(\bar u\) and \(\nabla\bar u\).
\end{proof}

For the frozen model equation, the Weiss and
Monneau energies obey exact identities; as we will see now, for our quasilinear equation, the
same identities persist up to an integrable error.

\begin{lem}\label{lem:degenerate-weiss-monneau}
Let \(\lambda=2m+2\), \(P=\bar c y_1^{2m}\), and define the Weiss and Monneau functionals
\[
W_\lambda(r,u):=
r^{-2\lambda}\int_{B_r}\bigl(|\nabla u|^2+2Pu\bigr)
-\lambda r^{-2\lambda-1}\int_{\partial B_r}u^2,\qquad M_\lambda(r,u,q):=r^{-2\lambda-1}\int_{\partial B_r}(u-q)^2,
\]
where \(q\ge0\) is a \(\lambda\)-homogeneous solution of  \(\Delta q=P\,\id_{\{q>0\}}\).
Then
\[
\frac{d}{dr}W_\lambda(r,\bar u)
=2r^{-2\lambda}\int_{\partial B_r}
\left(\partial_\nu\bar u-\frac{\lambda\bar u}{r}\right)^2\,d\Ha^1
+E_W(r),
\]
and
\[
\frac{d}{dr}M_\lambda(r,\bar u,q)
-\frac2r\bigl(W_\lambda(r,\bar u)-W_\lambda(1,q)\bigr)
=2r^{-2\lambda-1}\int_{B_r}
P\bigl(q\,\id_{\{\bar u=0\}}-\bar u\,\id_{\{q=0\}}\bigr)+E_M(r),
\]
with \(|E_W(r)|+|E_M(r)|\le C\).
In particular, if \(|\{q=0\}\cap B_1|=0\), then
\begin{equation*}
\frac{d}{dr}M_\lambda(r,\bar u,q)
\ge \frac2r\bigl(W_\lambda(r,\bar u)-W_\lambda(1,q)\bigr)+E_M(r).
\end{equation*}
\end{lem}

\begin{proof}
When \(v\) is a solution of the model equation
\[
\Delta v=P\,\id_{\{v>0\}},\qquad v\ge0,
\]
then the computation is the standard one of Weiss~\cite{Weiss1999} and
Monneau~\cite{Monneau2003}. Indeed, since \(P\) is homogeneous of degree
\(\lambda-2\), we have
\(y\cdot \nabla P=(\lambda-2)P\).
Differentiating \(W_\lambda\) and integrating by parts, all bulk terms cancel exactly due to this homogeneity, and the formula for \(\frac{d}{dr}W_\lambda(r,v)\) holds with \(E_W(r)=0\).

For Monneau's functional, put \(w=v-q\).  Since \(q\) is \(\lambda\)-homogeneous,
\[
\frac{d}{dr}M_\lambda(r,v,q)
=2r^{-2\lambda-1}\int_{\partial B_r}
w\left(\partial_\nu v-\frac{\lambda v}{r}\right),
\]
and the same integration-by-parts identity used in Weiss' formula gives
\[
\frac{d}{dr}M_\lambda(r,v,q)-\frac2r
\bigl(W_\lambda(r,v)-W_\lambda(1,q)\bigr)
=2r^{-2\lambda-1}\int_{B_r}
P\bigl(q\,\id_{\{v=0\}}-v\,\id_{\{q=0\}}\bigr).
\]
Now, by~\eqref{eq:ubar-pde}, the coefficient matrix in the second order part of the equation, namely \(A(\nabla \bar u + \nabla \varphi)\), equals the identity plus
\(O(|y|+|\nabla\bar u|)\), while (given the Lipschitz continuity of \(A\))
\[
|(A(\nabla \bar{u}+\nabla \varphi) - A(\nabla \varphi)) : D^2 \varphi| \leq C |\nabla \bar u| |D^2 \varphi| \leq C |y|^{4m+1}\, ,
\]
and moreover
\[
H(y_1)=P(y)+O(|y|^{2m+1}).
\]
Also, if we define \(u_r(y):=r^{-\lambda}\bar u(ry)\), then~\cref{lem:finite-order-upper-detachment} gives
\[
|u_r|+|\nabla u_r|+|D^2u_r|\le C,
\]
so the rescaled equation can be written as
\[
\Delta u_r
=P(y)\id_{\{u_r>0\}}+R_r(y),
\qquad |R_r(y)|\le Cr.
\]
The terms hidden in \(R_r\), including those coming from the nonconstant
principal matrix, enter the Rellich--Ne\v cas computation as additive bulk
errors bounded by \(Cr\).  Thus the usual proof, as in the perturbative
monotonicity formulas of~\cite{FocardiGelliSpadaro2015}, gives the estimates
from the statement.
\end{proof}

Next, we show that the positivity set cannot be too thin near a
free-boundary point, as it must carry a definite amount of weighted mass,
\(\int_{\{\bar u>0\}\cap B_r}|y_1|^{2m}\gtrsim r^{2m+2}\).

\begin{lem}\label{lem:weighted-density}
Set \(\lambda=2m+2\). There are \(\vartheta,r_*>0\) such that, for
\(0<r<r_*\),
\begin{equation*}
\sup_{\partial B_r}\bar u \geq \vartheta r^\lambda,\qquad
|\{\bar u>0\}\cap B_r|\ge \vartheta r^2,\qquad
\int_{\{\bar u>0\}\cap B_r}|y_1|^{2m}\,dy\ge \vartheta r^{\lambda}.
\end{equation*}
\end{lem}

\begin{proof}
We start with the lower bound on \(\sup_{\partial B_r}\bar u\) and we first
 recall~\cite[Theorem 1.1]{Yeressian2015}: if
\[
\Delta v=f\id_{\{v>0\}},\qquad
f(x)\ge\lambda_0 |(x_1,\ldots,x_j)|^\alpha,\qquad v\ge0,
\]
in some domain \(D\subset \mathbb R^n\),
then, for \(x^0\in\{v>0\}\) and \(B_\rho(x^0)\),
\begin{equation*}
\sup_{\{v>0\}\cap\partial B_\rho(x^0)}v
\ge v(x^0)+c\lambda_0 \rho^2
\bigl(\rho^\alpha+|(x^0_1,\ldots,x^0_j)|^\alpha\bigr).
\end{equation*}
Taking \(n=2\), \(j=1\) and \(\alpha=2m\), and letting \(x^0\to0\) through positivity
points, gives the lower bound
\begin{equation*}
\sup_{\partial B_r}v\ge c r^{2m+2}
\end{equation*}
for the model equation \(\Delta v=x_1^{2m}\id_{\{v>0\}}\) when
\(0\in\partial\{v>0\}\).

We now claim that the same proof as in~\cite{Yeressian2015} applies to~\eqref{eq:ubar-pde}.  Indeed, for a function \(w\), write
\[
\mathscr L w
:=A(\nabla\varphi+\nabla w):D^2w+
\bigl[A(\nabla\varphi+\nabla w)-A(\nabla\varphi)\bigr]:D^2\varphi .
\]
Fix \(y^0\in\{\bar u>0\}\), put \(a:=y^0_1\) and \(s:=y-y^0\),
and choose \(0<c_0<\bar c\). For
\(F(t):=c_0t^\lambda/[\lambda(\lambda-1)]\), adding sufficiently small
positive multiples of the harmonic polynomials
\[
|a|^{2m}(s_2^2-s_1^2),\qquad
(-1)^{m+1}\Re(s_1+i s_2)^\lambda
\]
to \(F(a+s_1)-F(a)-F'(a)s_1\) gives a polynomial \(\psi_0\) such that
\[
\Delta\psi_0=c_0y_1^{2m},\qquad
\psi_0(y^0)=0,\qquad
\psi_0(y)\ge c\bigl(|a|^{2m}|s|^2+|s|^\lambda\bigr).
\]
Set
\[
\psi(y):=\psi_0(y)-M\bigl(|a|^{2m+1}|s|^2+|s|^{\lambda+1}\bigr).
\]
Since \(\nabla\varphi(y)=O(|y_1|)\), \(D^2\varphi=O(1)\), and \(A(0)=I\),
the derivative bounds for \(\psi\) give
\begin{equation*}
|\mathscr L\psi-\Delta\psi|
\le C\bigl(|a|^{2m+1}+|s|^{2m+1}\bigr),
\end{equation*}
with \(C\) independent of \(M\), provided \(M(|a|+|s|)\) is sufficiently small.
Choosing \(M\) large and then shrinking the neighborhood, we obtain
\[
\mathscr L\psi\le c_0y_1^{2m}\le H(y_1),\qquad
\psi(y^0)=0,\qquad
\psi(y)\ge\frac c2\bigl(|a|^{2m}|s|^2+|s|^\lambda\bigr).
\]
The maximum principle applied to \(\bar u-\bar u(y^0)-\psi\) in
\(B_\rho(y^0)\cap\{\bar u>0\}\) therefore gives
\[
\sup_{\partial B_\rho(y^0)}\bar u
\ge\bar u(y^0)+c\rho^2\bigl(|y^0_1|^{2m}+\rho^{2m}\bigr).
\]
Letting \(y^0\to0\) through positivity points, we conclude that
\begin{equation}\label{eq:ubar-degenerate-nondegeneracy}
\sup_{\partial B_\rho}\bar u\ge c\rho^{2m+2}
\qquad\text{whenever }0<\rho<r_* .
\end{equation}
Next, we prove the third inequality. Assume by contradiction it fails. Then there is \(r_j\downarrow0\)
such that
\begin{equation}\label{e:ingtegral-vanishes}
r_j^{-2m-2}\int_{\{\bar u>0\}\cap B_{r_j}} |y_1|^{2m}\,dy\to0 .
\end{equation}
Set \(u_j(x):=r_j^{-\lambda}\bar u(r_jx)\).   Note that the sequence \(u_j\) is uniformly bounded in \(C^{1,1}(B_{1/2})\) (thanks to~\cref{lem:finite-order-upper-detachment}) and solves
\begin{align*}
&\; A \bigl(\nabla\varphi(r_jx)+r_j^{\lambda-1}\nabla u_j (x)\bigr):D^2u_j (x)
+  r_j^{2-\lambda}
\bigl[A \bigl(\nabla\varphi(r_jx)+r_j^{\lambda-1}\nabla u_j (x)\bigr)-A(\nabla\varphi(r_j x))\bigr]
:D^2\varphi(r_jx)\\
= &\; r_j^{2-\lambda} H(r_jx_1)\id_{\{u_j>0\}} (x)\, .
\end{align*}
We note that the coefficients converge uniformly to the identity, the gradient term is bounded by \(Cr_j\), and
\[
\int_{B_1}|r_j^{2-\lambda}H(r_jx_1)|\id_{\{u_j>0\}} (x)\,dx
= r_j^{-\lambda} \int_{B_{r_j}\cap \{\bar u >0\} } |H (y_1)| \,dy
\leq C r_j^{-2m-2} \int_{B_{r_j}\cap \{\bar u > 0\}} |y_1|^{2m}\, dy\; \stackrel{\eqref{e:ingtegral-vanishes}}{\to} 0 \, .
\]
Thus, passing to a subsequence, \(u_j\to w\) locally uniformly
with
\[
\Delta w=0,\qquad w\ge0,\qquad w(0)=0 .
\]
This implies \(w\equiv0\), and therefore
\[
r_j^{-\lambda}\sup_{B_{r_j/2}}\bar u=\sup_{B_{1/2}}u_j\to0,
\]
contradicting~\eqref{eq:ubar-degenerate-nondegeneracy}.

Having proven the third estimate the second one follows immediately from the elementary inequality
\[
\int_{\{\bar u>0\}\cap B_r}|y_1|^{2m}\,dy
\le r^{2m}|\{\bar u>0\}\cap B_r|\, . \qedhere
\]
\end{proof}

\begin{remark}\label{rem:polynomial-monneau}
We shall also use the polynomial Monneau comparison
\[
p_*(y):=\frac{\bar c\,y_1^{2m+2}}{(2m+2)(2m+1)},\qquad \Delta p_*=\bar c\,y_1^{2m}.
\]
The comparison has zero set of measure zero, so~\cref{lem:degenerate-weiss-monneau} gives a genuine Monneau monotonicity for
\(p_*\), up to an integrable perturbative error. Moreover, since the error \(E_W\) in the Weiss formula of~\cref{lem:degenerate-weiss-monneau} is integrable, \(W_\lambda(r,\bar u)\) is monotone up to an integrable error, and therefore  the
\(\lambda\)-density
\[
\Theta_\lambda(0):=\lim_{r\downarrow0}W_\lambda(r,\bar u)
\]
exists.
\end{remark}

Combining the upper bound with the nondegeneracy, we can show that the rescalings \(r^{-\lambda} \bar u(r\,\cdot)\) converge to a
\(\lambda\)-homogeneous solution of the frozen equation \(\Delta
q=y_1^{2m}\id_{\{q>0\}}\), with \(\partial_{y_2} q\ge0\).

\begin{prop}\label{prop:homogeneous-blowup}
Every sequence \(r_j\downarrow0\) has a subsequence such that
\[
\bar u_{r_j}(y):=\frac{\bar u(r_jy)}{r_j^{2m+2}}
\]
converges in \(C^1_{\mathrm{loc}}\) to \(\bar c\, q_{2m}\), where \(q_{2m}\) is homogeneous
of degree \(2m+2\) and satisfies
\[
\Delta q_{2m}=y_1^{2m}\id_{\{q_{2m}>0\}},\qquad
\partial_{y_2} q_{2m}\ge0.
\]
\end{prop}

\begin{proof}
As in the proof of~\cref{lem:weighted-density}
we have
\[
A(\nabla\varphi(r_jy)+r_j^{2m+1}\nabla\bar u_{r_j} (y)):D^2\bar u_{r_j} (y) +O(r_j)
=r_j^{-2m}H(r_jy_1)\id_{\{\bar u_{r_j}>0\}} (y)\, ,
\]
while the coefficients of the elliptic second order part of the equation converge to the frozen matrix, and
\[
r_j^{-2m}H(r_jy_1)=\bar c\,y_1^{2m}+O(r_j|y_1|^{2m+1}).
\]
Thus, thanks to~\cref{lem:finite-order-upper-detachment} and~\cref{lem:weighted-density}, every limit \(q\) is nonzero and solves
\[
\Delta q=\bar c\,y_1^{2m}\id_{\{q>0\}}.
\]
Also, applying Weiss' formula (see~\cref{lem:degenerate-weiss-monneau}) on the interval \(r \in (r_j\rho,r_jR)\), rescaling, and
letting \(j\to\infty\), gives
\[
\int_\rho^R s^{-2\lambda}\int_{\partial B_s}
\left(\partial_\nu q-\frac{\lambda q}{s}\right)^2\,d\Ha^1\,ds=0 .
\]
This implies that \(r\,\partial_r q=\lambda q\), from which we infer that \(q\) is homogeneous of degree
\(\lambda\).  Finally \(\partial_{y_2} \bar u \ge 0\) implies \(\partial_{y_2} \bar u_j \geq 0\), which, by
the estimates in~\cref{lem:obstacle}, passes to the limit.
\end{proof}

Observe now that any possible limit as \(r_j\downarrow 0\) of a sequence \(\bar u_j\) is, up to a constant, a homogeneous quadratic solution \(q\) of the classical obstacle problem \(\Delta q = \id_{\{q>0\}}\) satisfying the additional constraint that \(\partial_{y_2} q \geq 0\). Thanks to the classical classification theorem of Caffarelli, we then conclude the following.

\begin{cor}\label{lem:nondegenerate-points}
Assume \(m=0\) and let \(\bar{u}_{r_j}\), \(\bar{c}\), and \(q=q_0\) be as in~\cref{prop:homogeneous-blowup}. Then
either \(q (y) = \frac{\bar c}{2} (e\cdot y)_+^2\) for some \(e\in \mathbb S^1\) with \(e\cdot e_2 \geq 0\)  or \(q (y) = \frac{\bar c}{2} y_1^2\).
\end{cor}

Motivated by the previous result, we now classify (in our two-dimensional setting) all possible limits
when \(m\geq 1\). As we shall see, besides the whole-space and half-space polynomials, the only additional
possibility is a solution with positivity set which is a sector of opening \(\pi/(2m+1)\).

\begin{lem}\label{lem:sector-classification}
Assume \(m\ge1\).  The nontrivial homogeneous solutions of
\[
\Delta q=y_1^{2m}\id_{\{q>0\}},\qquad \partial_{y_2} q\ge0,
\]
are, up to multiplicative constants:
\begin{enumerate}
\item[(i)] the whole-space polynomial
\[
\frac{y_1^{2m+2}}{(2m+2)(2m+1)};
\]
\item[(ii)] one of the half-space polynomials
\[
\frac{((y_1)_+)^{2m+2}}{(2m+2)(2m+1)},\qquad
\frac{((y_1)_-)^{2m+2}}{(2m+2)(2m+1)}\, ,
\]
where \((y_1)_+ = \max \{y_1, 0\}\) and \((y_1)_- = - \min \{y_1, 0\}\);
\item[(iii)] a ``cone-type'' solution whose positivity set is a sector of opening \(\frac{\pi}{2m+1}\) centered on the positive semiaxis \(\{y_1= 0, y_2>0\}\).
\end{enumerate}
\end{lem}

\begin{proof}
The computation is the classical-obstacle analogue of the two-dimensional
homogeneous ODE classification in Appendix~A of~\cite{FocardiSpadaro2018}; see also~\cite{Yeressian2016}, where the homogeneous solutions are classified for the forcing term \(|x_1|\).
Let \(k=m+1\).

Consider a solution \(q\) and let \(S := \{q>0\}\). By the homogeneity of \(q\) this set is a union of sectors; we carry out the analysis on each connected component of \(S\), which, with a slight abuse of notation, we still denote by \(S\): it is either the whole plane or
a sector \(S= \{(r\cos \theta, r \sin \theta) : r>0\, , \theta_- < \theta < \theta_+\}\).
On \(S\) the solution \(q\) can be written as the sum of a homogeneous
harmonic polynomial and a specific homogeneous solution of the PDE:
\[
q=r^{2k}\phi(\theta),\qquad \phi(\theta)=\left(\frac{\cos^{2k}\theta}{2k(2k-1)}
-A\cos(2k\theta-\gamma)\right)\, ,
\]
where we can assume \(A\geq 0\) thanks to the arbitrariness of the phase \(\gamma\).

Differentiating in the vertical direction \(y_2\) the positivity condition
gives
\begin{equation}\label{e:monot-homog}
\partial_{y_2}q=2k A r^{2k-1}
\sin((2k-1)\theta-\gamma) \geq 0.
\end{equation}
We immediately conclude that, if the opening angle of the sector is larger than \(\frac{\pi}{2k-1}\), then \(A\) must vanish identically. In this case we see that we immediately fall into either (i) or (ii).

\medskip

Thus, from now on we can assume that \(S\) has opening angle at most \(\frac{\pi}{2k-1}\) and \(A> 0\).
In this case, we start by observing that on \(\{\theta=\theta_\pm\}\) we have \(\partial_{y_2} q =0\) and in particular
\begin{equation}\label{e:partialy_2_0}
\sin ((2k-1) \theta_\pm - \gamma) = 0\, ,
\end{equation}
so that~\eqref{e:monot-homog} gives that the opening angle of the sector is \emph{exactly} \(\frac{\pi}{2k-1}\):
\begin{equation}\label{e:half-period}
\theta_+-\theta_-=\frac{\pi}{2k-1}.
\end{equation}
The equations \(q=0\) on both rays then give
\begin{align}
\tfrac{1}{2k(2k-1)} \cos^{2k} (\theta_\pm)
 &= A \cos ((2k-1) \theta_{\pm} - \gamma) \cos \theta_\pm - A \sin ((2k-1) \theta_\pm -\gamma) \sin \theta_\pm
 \nonumber\\
 &= A \cos ((2k-1) \theta_{\pm} - \gamma) \cos \theta_\pm\, .\label{e:q=0}
\end{align}
Furthermore~\eqref{e:partialy_2_0} and~\eqref{e:half-period} imply that
\[
\cos ((2k-1) \theta_+-\gamma) = - \cos ((2k-1) \theta_-- \gamma) = \pm 1 .
\]
On the other hand, since we are assuming \(A>0\), the positivity of \(\sin ((2k-1) \theta -\gamma)\) in the range \(\theta \in (\theta_-, \theta_+)\) implies that \(\cos ((2k-1) \theta_+-\gamma) \leq \cos ((2k-1) \theta_--\gamma)\) hence
leaving only the possibility
\begin{equation*}
\cos ((2k-1) \theta_+-\gamma) = - \cos ((2k-1) \theta_-- \gamma) = - 1 .
\end{equation*}
In particular~\eqref{e:q=0} reduces to the system of equations
\begin{equation}\label{e:q=0-2}
\left\{
\begin{array}{l}
(\cos^{2k-1} \theta_+ + 2k (2k-1) A) \cos \theta_+ = 0\\
(\cos^{2k-1} \theta_- - 2k (2k-1) A) \cos \theta_- = 0
\end{array}
\right.
\end{equation}
We split the solutions of~\eqref{e:q=0-2} into three possibilities, depending on whether two, one, or none among \(\cos \theta_\pm\) vanish. The case in which both vanish is excluded, because it would imply \(\theta_+-\theta_-=\pi\), contradicting~\eqref{e:half-period} (recall that \(k=m+1\geq 2\)).

In the case in which none vanish, we have
\begin{equation}\label{e:cos-symmetry}
\cos^{2k-1}\theta_+=-\cos^{2k-1}\theta_- = - 2k (2k-1) A,
\end{equation}
Without loss of generality we can assume \(\theta_- \in [0, 2\pi)\), so that (taking into account~\eqref{e:half-period} and \(2k-1\geq 3\)) the possible solutions are
\[
\theta_\pm=\frac{\pi}{2}\pm\frac{\pi}{2(2k-1)}\,
\]
and
\[
\theta_\pm = \frac{3\pi}{2}\pm \frac{\pi}{2(2k-1)}\, .
\]
However by~\eqref{e:cos-symmetry} and (recalling that \(A>0\)) \(\cos \theta_+\) must be negative and so the second alternative is excluded. Note also that a sector of this form meets both half-planes \(\{y_1>0\}\) and \(\{y_1<0\}\); hence it is the only connected component of \(\{q>0\}\), since any other component would be either a second such sector or (when \(A=0\)) one of the half-planes \(\{y_1>0\}\), \(\{y_1<0\}\).

We have therefore found the solution described in point (iii) and we are thus only left with the task of excluding the case in which one among \(\cos \theta_\pm\) vanishes. This could be accomplished with further explicit computations; however, it is easier to just notice that the free boundary \(\partial \{q>0\}\) would consist of a vertical half-line and one half-line which is not vertical and thus the same argument employed in~\cref{lem:boundary-empty-side} (using Hopf's Lemma) applies then to our case as well, leading to a contradiction. The argument gives in fact a more general conclusion which, for convenience, we record in the next lemma.
\end{proof}

\begin{lem}\label{lem:nonvertical-cones}
Let \(m\ge 0\) and let \(q\ge0\) be a nontrivial homogeneous solution of \(\Delta q=y_1^{2m}\id_{\{q>0\}}\) with \(\partial_{y_2} q\ge0\) and \(\partial_{y_2} q \not\equiv0\). Then
no boundary ray of the positivity set \(\{q>0\}\) is vertical.
\end{lem}

\begin{proof}
Inside the positivity sector \(\partial_{y_2} q\) is harmonic because
differentiating \(\Delta q =y_1^{2m}\) in the \(y_2\) direction gives
\(\Delta(\partial_{y_2} q)=0\). Thus, by the strong maximum principle,
\(\partial_{y_2} q>0\) unless it is identically zero.

Suppose a boundary ray \(L\) were vertical and choose a point
\(p\in L\setminus\{0\}\).  Along a free-boundary ray one has
\(\nabla q =0\), therefore
\[
w:=\partial_{y_2} q =0,\qquad \partial_{y_1} q =0
\qquad\text{on }L .
\]
Since the tangent to \(L\) is \(e_2\), we get
\[
\partial_{y_2} (\partial_{y_1} q) (y) =0  \qquad \mbox{for \(y\in L\).}
\]
As in the argument of~\cref{lem:boundary-empty-side}, since \(q\) is only \(C^{1,1}\), we understand the above as differentiating the trace of \(\partial_{y_1} q\) from the side of the positivity sector.
On the other hand, since \(w\) is positive and harmonic in the sector,
Hopf's lemma  gives
\[
0<\partial_\nu w(p)=\pm\partial_{y_1} \partial_{y_2} q(p),
\]
where \(\nu=\pm e_1\) is the inward unit normal to the vertical side of the positivity sector. The two equations above provide the desired contradiction.
\end{proof}

Since the homogeneous models form a finite (hence isolated) family, the
blow-up is in fact unique (see~\cite{Yeressian2016} for a similar argument). We summarize this in the next result.

\begin{lem}\label{lem:polynomial-asymptotics}
For \(m\ge1\), let \(\mathscr H_m\) be the finite set of normalized
models in~\cref{lem:sector-classification}, and put
\[
\mathscr H_{m,\bar c}:=\{\bar c q:q\in\mathscr H_m\}.
\]
If there is a positive \(\bar c\) such that, along some sequence \(r_j\downarrow0\), the rescalings
\[
\bar u_{r_j}(x)=r_j^{-\lambda}\bar u(r_jx),\qquad \lambda=2m+2
\]
converge to \(q\in\mathscr H_{m,\bar c}\), then \(\bar u_r\to q\) in \(C^1_{\mathrm{loc}}\) as \(r\to 0\).
\end{lem}

\begin{proof}
By the compactness above and the homogeneous classification, every sequence
\(r_j\downarrow0\) has a subsequence for which \(\bar u_{r_j}\) converges to an
element of the finite set \(\mathscr H_{m,\bar c}\).  Thus the cluster set of
\(\bar u_r\) as \(r\downarrow0\) is contained in \(\mathscr H_{m,\bar c}\).  On the
other hand, the standard Monneau continuity argument (equivalently, the
continuity of the map \(r\mapsto \bar u_r\) in \(C^1_{\mathrm{loc}}\), together with
the compactness estimates) shows that this cluster set is connected.  Since
\(\mathscr H_{m,\bar c}\) is finite, the cluster set consists of a single element.
If one subsequence converges to \(q\), that element is \(q\), and the full
family converges to \(q\).
\end{proof}

The next estimate allows us to exclude the flat polynomial models: if
the positivity set contains a cone of opening close to \(\pi\), the vertical
derivative \(\chi=\partial_{y_2}\bar u\) decays no faster than \(t^{\beta}\) with \(\beta\)
near \(1\) along the axis. This is, in fact, a quantitative version of the argument used in the proof of~\cref{lem:boundary-empty-side}, which excluded that the set \(\partial \{\bar u > 0\}\)
can contain a vertical segment.

\begin{lem}\label{lem:cone-barrier-chi}
Suppose that for some angle
\(\theta\in(0,\pi)\), direction \(\sigma\in\{-1,1\}\), and radius \(r_0>0\),
the cone
\[
K_{\theta,r_0}^{\sigma}:=
\{r(\sigma\cos\phi,\sin\phi):0<r<r_0,\ |\phi|<\theta/2\}
\]
is contained in one connected component of \(\{\bar u>0\}\).  Then, after
possibly decreasing \(r_0\), there are \(c>0\) and
\(\beta=\pi/\theta+o_{r_0}(1)\) such that
\[
\partial_{y_2}\bar u(t\sigma e_1)\ge c\,t^\beta,\qquad 0<t<r_0/2.
\]
In particular, by taking \(\theta\) sufficiently close to \(\pi\), one can make
\(\beta<1+\varepsilon\).
\end{lem}

\begin{proof}
In the positivity component, by~\cref{prop:strict-monotonicity} we have that \(\chi:=\partial_{y_2}\bar u>0\)  solves a uniformly elliptic
linear equation
\[
L\chi:=a_{ij}(y)\partial_{ij}\chi+b_i(y)\partial_i \chi=0
\]
with \(a_{ij}(0)=\delta_{ij}\), \(a_{ij}\) continuous, and \(b_i\) bounded.
Fix \(\varepsilon>0\), choose \(\theta'<\theta\) such that
\(a:=\pi/\theta'<\beta:=\pi/\theta+\varepsilon\), and set
\[
d:=\beta^2-a^2,\qquad \eta:=\frac{d}{4a^2},\qquad
b(r,\phi):=r^\beta\bigl(\cos(a\phi)+\eta\cos^2(a\phi)\bigr).
\]
In \(K_{\theta',r_0}^{\sigma}\), the function \(b\) is positive, vanishes
on the sides, and satisfies \(\Delta b\ge(d/2)r^{\beta-2}\).
Since \(a_{ij}\to\delta_{ij}\) and \(b_i\) are bounded, decreasing \(r_0\)
gives \(Lb\ge(d/4)r^{\beta-2}>0\).
The closed outer arc of the smaller cone at radius \(r_0/2\) lies in
\(\{\bar u>0\}\), so \(\chi\) has a positive minimum there.
Hence a sufficiently small positive multiple of \(b\) is bounded above
by \(\chi\) on the boundary, and the maximum principle gives the same
inequality inside the cone. Along its axis this yields
\(\chi(t\sigma e_1)\ge c_\varepsilon t^{\pi/\theta+\varepsilon}\),
which proves the conclusion since \(\varepsilon>0\) is arbitrary.
\end{proof}

We can now show that the
whole-space and half-space polynomial alternatives of the degenerate
classification, and the corresponding quadratic polynomial alternatives when
\(m=0\), do not occur.

\begin{lem}\label{lem:exclude-polynomial}
Let \(m\geq 0\), \(\bar u\), \(\bar u_{r_j}\), \(\bar c\), and \(q_{2m}\) be as in~\cref{prop:homogeneous-blowup}. For \(m\ge 0\) the whole-space and half-space polynomial alternatives do not occur. More precisely:
\begin{itemize}
\item[(a)] If \(m=0\), then \(q_0 (y) = \tfrac{1}{2} \bar c (e\cdot y)_+^2\) for some \(e\in \mathbb S^1\) with \(e\cdot e_2 >0\).
\item[(b)] If \(m\geq 1\), then \(q:=\frac{q_{2m}}{\bar c}\) is a solution as in~\cref{lem:sector-classification}(iii).
\end{itemize}
\end{lem}

\begin{proof} Without loss of generality we assume \(\bar c =1\).
Let \(\chi:=\partial_{y_2}\bar u\). We first consider the case \(m\geq 1\). If one of the polynomial
alternatives occurred with \(m\ge1\) along a sequence \(r_j\downarrow0\), the
corresponding homogeneous limit supplied by~\cref{prop:homogeneous-blowup}
would be independent of \(y_2\). Also, by the uniform convergence
\(\bar u_r\to q\), given any \(\theta<\pi\) and \(\sigma\in\{-1,1\}\) corresponding to a positive
side of the polynomial model,
the polynomial profile is bounded below by a positive constant on
\[
S_\theta^\sigma:=\{(\sigma\cos\phi,\sin\phi):|\phi|\le\theta/2\}.
\]
Then~\cref{lem:cone-barrier-chi} gives
\[
\chi(t\sigma e_1)\ge c t^\beta,\qquad 0<t<r_0/2,
\]
with \(\beta<1+\varepsilon\). Choosing \(\varepsilon>0\) so small that
\(\beta<2m+1\), this  contradicts
\(\chi(t\sigma e_1)=O(t^{2m+1})\).

For \(m=0\), the inward curvature is strictly negative near the base point,
so Simon's boundary regularity theorem~\cite{Simon1976} gives local
Lipschitz continuity of the trace (see also~\cite[p.~98]{Lin87}). Hence the
contact set in the obstacle chart contains a cone of positive opening
with axis \(-e_2\), and every blow-up vanishes on this cone.
By~\cref{lem:nondegenerate-points}, this excludes the whole-space quadratic
model and the vertical half-plane profiles, proving (a).
\end{proof}

Finally, Caffarelli's regular
free-boundary theory~\cite{Caffarelli1977,Caffarelli1998,
PetrosyanShahgholianUraltseva2012} and the analytic bootstrap from~\cite{KinderlehrerNirenberg1977,KinderlehrerNirenbergSpruck1978} yield analyticity of
the free boundary near the nondegenerate points where \(m=0\), since we have shown that the only possible blow-ups are the half-plane solutions
\(\frac{\bar{c}}{2} (y\cdot e)_+^2\) with \(e\cdot e_2>0\). For the trace itself, this case (\(\kappa<0\) and \(h\neq g\)) is also contained in~\cite[Theorem~1 and Corollary~1]{Lin87}.

\begin{cor}\label{lem:regular-points}
If \(m=0\) then the corresponding
free-boundary point in \(\partial \{\bar u >0\}\) is regular and the free boundary \(\partial \{\bar u =0\}\) is
\(C^\omega\) in a neighborhood.
\end{cor}

\section{Boundary blow-ups}\label{sec:boundary-blowup}

We now study free-boundary points lying on the fixed boundary
\(\partial_b \Lambda\), namely the endpoints of the detachment set, where
\(\tilde h\) meets \(\tilde g\).  We keep the setup of~\cref{sec:obstacle,sec:sign-reduction}: \(\bar u\ge0\) solves~\eqref{eq:ubar-pde} on \(\Lambda=\{y_2>\tilde g(y_1)\}\), it vanishes on
\(\partial_b \Lambda=\{y_2=\tilde g(y_1)\}\), satisfies \(\partial_{y_2}\bar u\ge0\),
and its free boundary is the graph \(\{y_2=\tilde h(y_1)\}\) of the trace.  We
write the first nonzero term of the analytic forcing as
\[
H(y_1)=c\,y_1^\ell+O(y_1^{\ell+1}),\qquad c\ne0.
\]

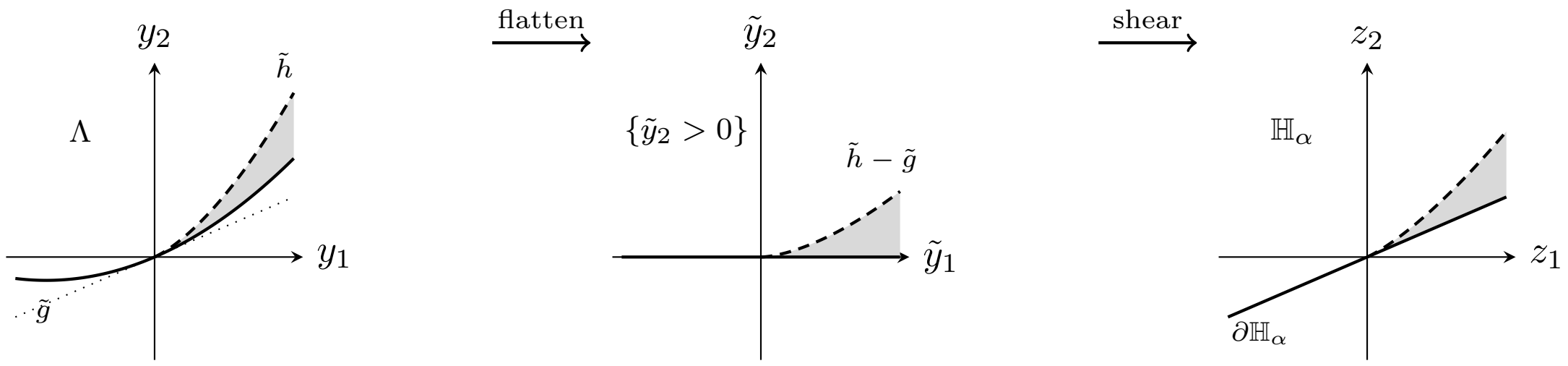
\begin{figure}[hbtp]
\centering
\begin{tikzpicture}
  \pgfplotsset{
    chartpanel/.style={
      width=4.3cm, height=4.3cm, scale only axis=false,
      axis lines=middle, xtick=\empty, ytick=\empty,
      xmin=-1.6, xmax=1.6, ymin=-0.9, ymax=1.7,
      samples=80, clip=false, anchor=south west,
      label style={font=\footnotesize},
      every axis x label/.style={at={(current axis.right of origin)}, anchor=west},
      every axis y label/.style={at={(current axis.above origin)}, anchor=south},
    }
  }
  \begin{axis}[chartpanel, at={(0cm,0cm)}, xlabel={\(y_1\)}, ylabel={\(y_2\)}]
    \addplot[name path=G1, thick, domain=-1.5:1.5] {0.35*x+0.15*x^2};
    \addplot[name path=H1, thick, densely dashed, domain=0:1.5] {0.35*x+0.15*x^2+0.3*x^1.6};
    \addplot[black!15] fill between[of=H1 and G1, soft clip={domain=0:1.5}];
    \addplot[thin, dotted, domain=-1.5:1.5] {0.35*x};
    \node[font=\footnotesize] at (axis cs:-0.8,1.1) {\(\Lambda\)};
    \node[font=\scriptsize, anchor=north] at (axis cs:-1.2,-0.25) {\(\tilde g\)};
    \node[font=\scriptsize, anchor=south] at (axis cs:1.4,1.45) {\(\tilde h\)};
  \end{axis}
  \begin{axis}[chartpanel, at={(5.55cm,0cm)}, xlabel={\(\tilde y_1\)}, ylabel={\(\tilde y_2\)}]
    \addplot[name path=G2, thick, domain=-1.5:1.5] {0};
    \addplot[name path=H2, thick, densely dashed, domain=0:1.5] {0.3*x^1.6};
    \addplot[black!15] fill between[of=H2 and G2, soft clip={domain=0:1.5}];
    \node[font=\footnotesize] at (axis cs:-0.8,1.1) {\(\{\tilde y_2>0\}\)};
    \node[font=\scriptsize, anchor=south] at (axis cs:1.3,0.62) {\(\tilde h-\tilde g\)};
  \end{axis}
  \begin{axis}[chartpanel, at={(11.1cm,0cm)}, xlabel={\(z_1\)}, ylabel={\(z_2\)}]
    \addplot[name path=G3, thick, domain=-1.5:1.5] {0.35*x};
    \addplot[name path=H3, thick, densely dashed, domain=0:1.5] {0.35*x+0.3*x^1.6};
    \addplot[black!15] fill between[of=H3 and G3, soft clip={domain=0:1.5}];
    \node[font=\footnotesize] at (axis cs:-0.8,1.1) {\(\mathbb H_\alpha\)};
    \node[font=\scriptsize, anchor=north] at (axis cs:-1.15,-0.45) {\(\partial\mathbb H_\alpha\)};
  \end{axis}
  \draw[->, thick] (4.45cm,2.9cm) -- node[above, font=\scriptsize]{flatten} (5.35cm,2.9cm);
  \draw[->, thick] (10.0cm,2.9cm) -- node[above, font=\scriptsize]{shear} (10.9cm,2.9cm);
\end{tikzpicture}
\caption{The change of variables at a boundary contact point: the flattening \(\tilde y_2=y_2-\tilde g(y_1)\), then the shear \(z_2=\tilde y_2+a\tilde y_1\), with \(a=\tilde g'(0)=\tan\alpha\). Solid: the fixed boundary; dashed: the free boundary; shaded: the contact set.}\label{f:boundary-chart}
\end{figure}

Fix one such point and translate it to the origin.  To flatten the fixed
boundary we use the analytic change of variables
\[
\tilde y_1=y_1,\qquad \tilde y_2=y_2-\tilde g(y_1),
\]
which transforms \(\partial_b \Lambda\) to \(\{\tilde y_2=0\}\) and \(\Lambda\) to the
half-plane \(\{\tilde y_2>0\}\). Note that, since it preserves the vertical direction,
\(\partial_{\tilde y_2}\bar u\ge0\) and the free boundary remains a graph over
\(\tilde y_1\).  Put \(a:=\tilde g'(0)\).  Since \(\nabla\varphi(0)=0\), the
principal part of~\eqref{eq:ubar-pde} at the origin is the Laplacian, which in
the \(\tilde y\)-variables becomes
\[
L_0=(\partial_{\tilde y_1}-a\partial_{\tilde y_2})^2+
\partial_{\tilde y_2}^2 .
\]
The linear map
\[
z_1=\tilde y_1,\qquad z_2=\tilde y_2+a\tilde y_1
\]
reduces \(L_0\) to the Laplacian, preserves the vertical derivative
\(\partial_{\tilde y_2}=\partial_{z_2}\), and keeps the forcing variable
\(\tilde y_1=z_1\).  Also, this map sends \(\{\tilde y_2>0\}\) to the half-plane
\[
\mathbb H_\alpha=\{\theta\in(\alpha,\alpha+\pi)\},
\]
where
\[
\alpha:=\arctan a\in(-\pi/2,\pi/2).
\]
In particular \(\cos\alpha>0\), cf.~\cref{f:boundary-chart}.  We keep the notation \(\bar u\) for the
transformed function.

\begin{lem}\label{lem:one-sided-positive-sign}
Assume that \(0\in\partial_b \Lambda\) is a boundary free-boundary point and that
detachment occurs on one side of \(0\) along \(\partial_b \Lambda\).  Then the
first nonzero term of \(H\) has the sign for which \(H>0\) on that side.
\end{lem}

\begin{proof}
Write \(H(y_1)=c\,y_1^\ell+O(y_1^{\ell+1})\).  If the leading term were
negative on the detaching side, then, after shrinking the interval on that
side, \(H<0\) there.  This would mean that the inward curvature is
nonnegative, so~\cref{prop:attachment} gives \(\tilde h=\tilde g\) there,
contradicting detachment.
\end{proof}

The boundary blow-ups are the half-plane analog of the interior ones of~\cref{prop:homogeneous-blowup}: the rescalings converge to a homogeneous
solution vanishing on the flat Dirichlet side.

\begin{lem}\label{lem:boundary-estimates}
At a boundary free-boundary point, with \(d=\ell+2\) and with the sign fixed by~\cref{lem:one-sided-positive-sign}, there are constants \(c_0,C,r_0>0\) such
that, for \(0<r<r_0\),
\[
c_0r^d\le\sup_{B_r\cap\mathbb H_\alpha}\bar u\le Cr^d,\qquad
\sup_{B_r\cap\mathbb H_\alpha}|\nabla\bar u|\le Cr^{d-1},\qquad
\sup_{B_r\cap \mathbb H_\alpha}|D^2 \bar u| \le C r^{d-2}\, .
\]
Moreover the corresponding half-ball Weiss functional of homogeneity \(d\) is
quasi-monotone, with the same type of integrable error as in~\cref{lem:degenerate-weiss-monneau}.
\end{lem}

\begin{proof}
First of all we recall that \(\mathscr L\) has divergence form. Indeed, for a \(C^2\) function \(f\) the mean curvature of its graph is \(H_f=\operatorname{div}\bigl(D_pF(\nabla f)\bigr)\), with \(F\) as in~\eqref{e:variational-ineq}, and therefore
\[
\mathscr L w=H_{\varphi+w}-H_\varphi=\operatorname{div}\bigl(D_pF(\nabla\varphi+\nabla w)-D_pF(\nabla\varphi)\bigr).
\]
Denote by \(y=Y(z)\) the inverse of the change of variables of this section, namely the flattening \(\tilde y\) followed by the linear map \(z\). Note that \(Y\) is an analytic diffeomorphism with \(\det DY=1\) (since both maps are shear), so in the variables \(z\) we still have
\[
 \mathscr L w=\operatorname{div}\bigl(\Phi(z,\nabla w)-\Phi(z,0)\bigr),
\]
where \(\Phi\) is analytic and \(D_\xi\Phi\) is symmetric and uniformly elliptic on bounded sets. Consequently, for any two functions \(w_1,w_2\) of class \(C^1\),
\begin{equation}\label{eq:linearized-difference}
\mathscr L w_1-\mathscr L w_2=\operatorname{div}\bigl(A\,\nabla(w_1-w_2)\bigr),\qquad A:=\int_0^1D_\xi\Phi\bigl(z,\nabla w_2+s\nabla(w_1-w_2)\bigr)\,ds .
\end{equation}
We now prove the upper bound in the first estimate. Fix \(r\in(0,r_0/2)\), set \(\mathcal B:=B_{2r}\cap\mathbb H_\alpha\) and \(\delta(z):=\dist(z,\partial\mathbb H_\alpha)\), and let \(e\) be the inner unit normal to \(\partial\mathbb H_\alpha\). Since \(|H|\le C|z|^\ell\), we have \(|\mathscr L\bar u|\le Cr^{d-2}\) in \(\mathcal B\). Let \(\hat u\) be the solution of \(\mathscr L\hat u=0\) in \(\mathcal B\) with \(\hat u=\bar u\) on \(\partial\mathcal B\). Since \(\bar u\in C^{1,1}\), classical estimates for quasilinear equations in divergence form imply
\[
 \|\nabla\hat u\|_{L^\infty(B_{3r/2}\cap\mathbb H_\alpha)}+r^{\alpha}[\nabla\hat u]_{C^{\alpha}(B_{3r/2}\cap\overline{\mathbb H_\alpha})}\le C,
\]
 with \(C\) independent of \(r\). By~\eqref{eq:linearized-difference}, \(\hat u\) solves the linear equation \(\operatorname{div}(A_0\nabla\hat u)=0\), where \(A_0\) is the matrix \(A\) of~\eqref{eq:linearized-difference} with \(w_1=\hat u\) and \(w_2=0\); and \(w:=\bar u-\hat u\) solves
\[
\operatorname{div}(A\nabla w)=\mathscr L\bar u\quad\text{in }\mathcal B,\qquad w=0\quad\text{on }\partial\mathcal B.
\]
 The matrices \(A_0\) and \(A\) are uniformly elliptic and, after the rescaling \(z\mapsto rz\), their \(C^\alpha\) norms on \(B_{3/2}\cap\overline{\mathbb H_\alpha}\) are bounded independently of \(r\).
Thus, \(L^\infty\) estimate for divergence-form equations gives
 \[
  \sup_{\mathcal B}|w|\le Cr^2\sup_{\mathcal B}|\mathscr L\bar u|\le Cr^d,
 \]
 and the \(C^{1,\alpha}\) estimates up to the flat part of \(\partial\mathcal B\) give
 \[
  |\nabla w|\le C\bigl(r^{-1}\sup_{\mathcal B}|w|+r\sup_{\mathcal B}|\mathscr L\bar u|\bigr)\le Cr^{d-1}
 \]
 in \(B_r\cap\mathbb H_\alpha\). Since \(w=0\) on \(\partial\mathbb H_\alpha\), we get \(|w|\le Cr^{d-1}\delta\) in \(B_r\cap\mathbb H_\alpha\).

By the maximum principle, \(\hat u\ge0\). Now, set \(M:=\hat u(re)\). The Harnack inequality and the Carleson estimate~\cite[Theorem 1.1]{CaffFabeMortSals81} give \(\hat u\le CM\) in \(B_{3r/2}\cap\mathbb H_\alpha\); thus, \(C^{1,\alpha}\) estimates up to the flat boundary (where \(\hat u=0\))  give
 \[
  \hat u\le C(M/r)\,\delta  \qquad \text{in \( B_r\cap\mathbb H_\alpha\).}
 \]
 Conversely, the Hopf lemma (which holds for divergence-form operators with H\"older continuous coefficients) together with the Harnack inequality yields \(\hat u\ge c\,(M/r)\,\delta\) in \(B_r\cap\mathbb H_\alpha\).

Since by assumption there is \(z_0\in B_r\cap\mathbb H_\alpha\) with \(\bar u(z_0)=0\), we get \(\hat u(z_0)=-w(z_0)\le Cr^{d-1}\delta(z_0)\), so the lower bound for \(\hat u\) at \(z_0\) gives \(M\le Cr^d\). Therefore, in \(B_r\cap\mathbb H_\alpha\),
\[
\bar u=\hat u+w\le C\,\frac Mr\,\delta+Cr^{d-1}\delta\le Cr^d .
\]

Having established the upper bound in the first estimate, the second upper bounds follow from classical elliptic theory. To get the third we can apply~\cite[Theorem 4.1]{FocardiGeraciSpadaro2017} to the function \(r^{-d} \bar{u} (r z)\).

For the lower bound in the first estimate, we distinguish the parity of \(\ell\). If \(\ell\) is even, then \(H\ge0\) near \(0\) by~\cref{lem:one-sided-positive-sign}, and the argument used in~\cref{lem:weighted-density} applies in half-balls, again with barriers
vanishing on the flat side.

If \(\ell\) is odd, then \(H\) changes sign at \(0\). By~\cref{lem:one-sided-positive-sign}, \(H>0\) on the detaching side, so \(H<0\) on the other side, where~\cref{prop:attachment} gives \(\tilde h=\tilde g\). Let \(\sigma\in\{-1,1\}\) be such that \(\sigma z_1>0\) on this attached side, and set \(D_r:=B_r\cap\mathbb H_\alpha\cap\{\sigma z_1> r\cos\alpha/4\}\). Then \(D_r\subset\{\bar u>0\}\) and \(H\le-c\,r^{\ell}\) on \(D_r\) for \(r\) small. By~\eqref{eq:linearized-difference} with \(w_1=\bar u\) and \(w_2=0\), we have
\[
\operatorname{div}(\bar a\nabla\bar u)=\mathscr L\bar u=H\qquad \text{in } D_r,
\]
with \(\bar a\) uniformly elliptic. Let \(\varpi_r\) solve \(\operatorname{div}(\bar a\nabla\varpi_r)=-1\) in \(D_r\) with \(\varpi_r=0\) on \(\partial D_r\); by scaling and uniform ellipticity, \(\sup_{D_r}\varpi_r\ge c\,r^2\). The function \(\bar u-c\,r^\ell\varpi_r\) satisfies \(\operatorname{div}\bigl(\bar a\nabla(\bar u-c\,r^\ell\varpi_r)\bigr)=H+c\,r^\ell\le0\) in \(D_r\) and is nonnegative on \(\partial D_r\), so the maximum principle gives \(\bar u\ge c\,r^\ell\varpi_r\) in \(D_r\), and therefore \(\sup_{B_r\cap\mathbb H_\alpha}\bar u\ge c\,r^{\ell+2}=c\,r^d\).

Finally, the Weiss computation can be performed in \(B_r\cap\mathbb H_\alpha\).
The contribution from the flat boundary vanishes because \(\bar u=0\) there,
and the errors coming from the analytic flattening, from the coefficients, and
from \(H-cz_1^\ell\) are \(O(1)\), as in~\cref{lem:degenerate-weiss-monneau}.
\end{proof}

\begin{lem}\label{lem:boundary-blowups}
At a boundary free-boundary point, every sequence \(r_j\downarrow0\) has a
subsequence such that
\[
u_{r_j}(z):=\frac{\bar u(r_jz)}{r_j^d}
\]
converges in \(C^1_{\mathrm{loc}}(\overline{\mathbb H_\alpha})\) to a nonzero homogeneous
function \(Q\) of degree \(d\) satisfying
\begin{equation}
\label{eq:Q}
\Delta Q=c\,z_1^\ell\id_{\{Q>0\}}\ \text{ in }\mathbb H_\alpha,\qquad
Q=0\ \text{ on }\partial\mathbb H_\alpha,\qquad \partial_{z_2}Q\ge0 .
\end{equation}
\end{lem}

\begin{proof}
By~\cref{lem:boundary-estimates}, the rescalings are uniformly bounded in
\(C^{1,1}_{\mathrm{loc}}(\overline{\mathbb H_\alpha})\).  The
rescaled equations have coefficients converging to the identity, while the
right-hand sides converge to
\[
c\,z_1^\ell\id_{\{Q>0\}} .
\]
Thus every limit satisfies the ``frozen'' equation with flat Dirichlet
condition.  The half-ball Weiss
quasi-monotonicity from~\cref{lem:boundary-estimates}, applied on
\((r_j\rho,r_jR)\) and then rescaled, gives homogeneity of degree \(d\), exactly
as in~\cref{prop:homogeneous-blowup}.  The inequality
\(\partial_{z_2}Q\ge0\) follows from~\cref{lem:obstacle}.
\end{proof}

The half-plane models are classified by the same angular ODE, now with one
endpoint pinned on the fixed Dirichlet boundary.  As in the interior, no
genuine free-boundary ray is vertical, with a quantitative angle gap.

\begin{lem}\label{lem:half-plane-classification}
Let \(d=\ell+2\) and \(Q=r^d\psi(\theta)\) be a nonzero \(d\)-homogeneous solution  of~\eqref{eq:Q}.
On each positivity sector \(S=\{(r\cos \theta, r\sin \theta) : \theta_- < \theta < \theta_+, r>0\}\subset \mathbb H_\alpha\) we have
\[
Q=r^d\left(c\,\frac{\cos^d\theta}{d(d-1)}
+A\cos(d\theta)+B\sin(d\theta)\right).
\]
If \(\theta_- \in \mathbb H_\alpha\) (or analogously \(\theta_+\)), we have
\(\psi(\theta_-)=\psi'(\theta_-)=0\); instead, if
\(\theta_- \in \partial\mathbb H_\alpha\), then \(\psi(\theta_-)=0\) (and analogously \(\theta_+\)).
Also, \(\theta_-,\theta_+ \neq \pi/2\), namely, no free-boundary ray
is vertical.
\end{lem}

\begin{proof}
The PDE gives
\[
\psi''+d^2\psi=c\cos^\ell\theta,\qquad \theta \in (\theta_-,\theta_+).
\]
If both endpoints are free, then
\(w:=\partial_{z_2}Q\) is positive and harmonic in the sector, unless it
vanishes identically. Thus, the Hopf argument of~\cref{lem:nonvertical-cones} then
excludes a vertical  free ray.

Assume now that one
endpoint is the fixed boundary ray \(\alpha\) and the other endpoint \(\beta\) is free.
Then, since
\[
(\cos^d\theta)''+d^2\cos^d\theta
=d(d-1)\cos^{d-2}\theta=d(d-1)\cos^\ell\theta,
\]
the ODE together with the conditions \(\psi(\beta)=\psi'(\beta)=0\) gives
\[
\psi(\theta)
=\frac{c}{d(d-1)}\Bigl[
\cos^d\theta-\cos^{d-1}\beta\,
\cos\bigl(d\theta-(d-1)\beta\bigr)\Bigr].
\]
Using that \(\psi(\alpha)=0\),
we obtain
\begin{equation}\label{eq:one-sided-angular-equation}
\cos^{\ell+1}\beta\,
\cos\bigl((\ell+2)\alpha-(\ell+1)\beta\bigr)
=\cos^{\ell+2}\alpha.
\end{equation}
By our choice of \(\alpha\), \(\cos\alpha>0\).  Substituting \(\beta=\pi/2\) would make the left-hand side
zero and the right-hand side nonzero, so a free-boundary ray cannot be
vertical.
\end{proof}

Note that for fixed \(\alpha \in (-\pi/2,\pi/2)\) the set of possible solutions to \eqref{eq:one-sided-angular-equation} is finite and depends only on \(\alpha\), in particular for any given \(\alpha\) the set of  \(\psi\) for which \(Q=r^{d}\psi(\theta)\) is an admissible blow-up is discrete.

\begin{lem}\label{lem:uniqueness-of-boundary-blow-ups}
 Assume that \(x_0\in\{g=h\}\), let \(\ell\ge 1\) and \(d=\ell+2\). Let \(\mathcal{K}\) be the set of limit points of \(r_j^{-d}\bar u(r_j\cdot)\) for subsequences \(r_j\to 0\). Then \(\mathcal{K}\) is a singleton.
\end{lem}

\begin{proof}
Since the map \((0,1)\ni r \mapsto r^{-d}\bar u(r\cdot)\) is continuous, \(\mathcal{K}\) is connected. Every \(Q\in\mathcal K\) has a genuine free ray. Indeed, if \(Q=r^d\psi(\theta)\) were positive on all of \(\mathbb H_\alpha\), then \(\psi=\frac{c}{d(d-1)}\cos^d\theta+R\cos(d\theta-\gamma)\) there, and since \(\partial_{z_2}(z_1^d)=0\) we would have \(\partial_{z_2}Q=-dR\,r^{d-1}\sin\bigl((d-1)\theta-\gamma\bigr)\), which is nonnegative on an interval of length \(\pi\) only if \(R=0\), because \((d-1)\pi\ge2\pi\); but then \(\psi=\frac{c}{d(d-1)}\cos^d\theta\) does not vanish on \(\partial\mathbb H_\alpha\), since \(\cos\alpha>0\). On the other hand, the above discussion says that \(\mathcal{K}\) is a finite discrete set that depends only on \(\alpha\) and \(c\) (and thus on  \(g\) and \(H\)). This immediately implies the conclusion.
 \end{proof}

\section{Proof of the main theorem}\label{sec:main-proof}

First of all, in the setup of~\cref{sec:obstacle} (in particular, the free boundary is the graph \(\{y_2=\tilde h(y_1)\}\)), we can apply the interior and boundary
classifications to obtain a uniform Lipschitz estimate.

\begin{prop}\label{prop:lipschitz}
The trace \(\tilde h\) is Lipschitz and admits one-sided
derivatives.
\end{prop}
\begin{proof}
\smallskip\noindent
\emph{Step 1: continuity.}
Since the contact set is relatively closed, \(\tilde h\) is
upper semicontinuous, so a discontinuity would produce a
nontrivial vertical segment of free boundary strictly above
the fixed boundary. Blowing up at an interior point of this
segment would give \(Q(e_2)=0\), contradicting
\cref{lem:exclude-polynomial,lem:sector-classification},
according to which every admissible interior model satisfies
\(Q(e_2)>0\). Hence \(\tilde h\) is continuous.

\smallskip\noindent
\emph{Step 2: Lipschitz continuity.}
Fix a point of the graph of \(\tilde h\), translate it to the
origin, and use the coordinates of~\cref{sec:boundary-blowup}
when it lies on the fixed boundary. If \(\tilde h=\tilde g\)
near this point, the conclusion follows from the regularity
of \(\tilde g\); otherwise, we claim that for some \(C,\rho>0\),
\[
|\partial_{y_1}\bar u|\le C\partial_{y_2}\bar u
\qquad\text{in }B_\rho\cap\{\bar u>0\}.
\]
Suppose by contradiction that there are \(y_j\to0\) in
\(\{\bar u>0\}\) such that
\[
|\partial_{y_1}\bar u(y_j)|>j\,\partial_{y_2}\bar u(y_j) \qquad \mbox{and} \qquad \partial_{y_2}\bar u(y_j)>0.
\]
Set \(r_j:=|y_j|\). By~\cref{prop:homogeneous-blowup,lem:boundary-blowups},
along a subsequence we have
\[
u_j(y):=r_j^{-\lambda}\bar u(r_jy)\longrightarrow Q,
\qquad z_j:=y_j/r_j\longrightarrow z\in\partial B_1,
\]
in \(C^1_{\mathrm{loc}}\), up to the fixed boundary when present,
with the appropriate homogeneity \(\lambda\).
The gradient inequality implies \(\partial_{y_2}Q(z)=0\).
Since \(\partial_{y_2}Q>0\) in each positivity sector, the
nondegeneracy barriers of~\cref{lem:weighted-density,lem:boundary-estimates},
applied in balls centered at \(z_j\), show that \(z\) lies on
a boundary ray of \(\{Q>0\}\), possibly on the fixed side.
By~\cref{lem:exclude-polynomial,lem:sector-classification,lem:half-plane-classification},
this ray is nonvertical; moreover, \(Q(z)=|\nabla Q(z)|=0\).
The angular equation then shows that the limiting forcing
\(f\) is positive at \(z\) and that \(Q\) has positive quadratic
growth along the ray, with inward normal \(\nu\) satisfying \(\nu_2>0\).

Let \(L_j\) be the linearization of the rescaled divergence-form
operator at \(u_j\). Since the flux and forcing depend only
on \(y_1\), the \(C^2\) convergence of the fluxes to the identity
and the uniform \(C^{1,1}\) estimates give
\[
L_j(\partial_{y_2}u_j)=0,\qquad
|L_j(\partial_{y_1}u_j)|\le C_0,\qquad
L_ju_j\ge c_0>0
\]
in \(B_\delta(z)\cap\{u_j>0\}\), for \(\delta>0\)
sufficiently small, and the coefficients of \(L_j\) are
uniformly elliptic and Lipschitz.
Choose \(K>(C_0+1)/c_0\). By the quadratic expansion of
\(Q\) along the ray, taking \(C\) sufficiently large
and decreasing \(\delta\) gives
\[
w_j^\pm:=C\partial_{y_2}u_j
\pm\partial_{y_1}u_j-Ku_j\ge-\varepsilon_j,
\qquad \varepsilon_j\to0,\qquad L_jw_j^\pm\le-1.
\]
On the free boundary \(w_j^\pm=0\), while on a fixed side
\(\{y_2=a y_1\}\) we have
\(w_j^\pm=(C\mp a)\partial_{y_2}u_j\ge0\), provided \(C>|a|\).
If \(w_j^\pm(x_0)<0\) at some
\(x_0\in B_{\delta/2}(z)\cap\{u_j>0\}\), then, for a small fixed
\(\beta>0\), the function
\[
w_j^\pm(y)+\beta|y-x_0|^2
\]
is strictly \(L_j\)-superharmonic and nonnegative on the
boundary of \(B_{\delta/2}(x_0)\cap\{u_j>0\}\) for large \(j\),
contradicting the minimum principle. Hence
\(|\partial_{y_1}u_j(z_j)|\le C\partial_{y_2}u_j(z_j)\),
which proves the claim.

Returning to the original coordinates and increasing \(C\)
if necessary, the directions \((1,C)\) and \((-1,C)\) preserve
\(\Lambda\). The derivative estimate gives monotonicity in
these directions, so the boundary of the epigraph
\(\{\bar u>0\}=\{y_2>\tilde h(y_1)\}\) is a Lipschitz graph.

\smallskip\noindent
\emph{Step 3: one-sided derivatives.}
At points where \(\tilde h=\tilde g\) in a neighborhood,
the conclusion follows from the regularity of \(\tilde g\).
Otherwise the blow-up is unique by
\cref{lem:polynomial-asymptotics,lem:regular-points,lem:uniqueness-of-boundary-blow-ups},
with the nondegenerate boundary case given by the quadratic
classification.
By the Lipschitz bound, the rescaled graphs are locally
compact in the uniform topology, and every subsequential
limit coincides with the lower boundary of \(\{Q>0\}\):
indeed, \(Q\) vanishes below the limiting graph and
satisfies \(\Delta Q=f\) above it, where \(f\neq0\)
away from \(\{y_1=0\}\).
Hence the graphs converge to the two boundary rays of
a single model, including fixed-boundary rays, whose
slopes give the one-sided derivatives of \(\tilde h\).
\end{proof}

We can now begin proving~\cref{thm:main}. We do this in steps.

\begin{proof}[Proof of~\cref{thm:main}(i)--(ii)]
At points of \(\Gamma_r\), \(h=g\).  At the other points \(x\), we have \((x,g(x))\in\overline{\Gamma^+\cup\Gamma^-}\), and~\cref{lem:obstacle} represents \(h\) as the free boundary \(\tilde h\) of an obstacle problem, and by~\cref{prop:lipschitz} this function is locally Lipschitz and admits one-sided derivatives. Hence, by compactness,
\(h\) is globally Lipschitz and admits one-sided derivatives at every point.

From this, the relation between the support of the measure \(\Ha^2 \res (\partial^*E\cap(\partial\Omega\times\R))\) and the graphs of \(h\) and \(g\) is immediate.
This proves (i).

We next address the first part of (ii). Fix \(p\in\{h=g\}\setminus\Sigma\) and assume that detachment set $\{h\neq g\}$ accumulates
at \(p\) (otherwise \(h=g\) near \(p\) and \(h\) is analytic there, hence there is nothing to prove).
Then \(\kappa(p)<0\) by~\cref{prop:attachment}, so in the
coordinates of~\cref{sec:boundary-blowup}, every blow-up has the form
\(Q=r^2\psi(\theta)\) with constant forcing \(c=H(0)>0\). We first claim that both endpoints are fixed.
Indeed, at a free endpoint \(\beta\in(\alpha,\alpha+\pi)\) the angular equation
and zero Cauchy data of~\cref{lem:half-plane-classification} would give
\(\psi(\theta)=\tfrac c2\sin^2(\theta-\beta)\), whose consecutive zeros are
\(\pi\) apart, forcing the sector to coincide with \(\mathbb H_\alpha\) and
\(\beta\) to be fixed. Thus \(Q\) is a nonnegative quadratic polynomial vanishing
on the fixed line, so \(\Delta Q=c\) gives \(Q(z)=\tfrac c2(z\cdot n)^2\),
where \(n\) is the inward unit normal to \(\mathbb H_\alpha\).

Now, after translating \(p\) to the origin, we denote by \(\tilde h(s)\) and
\(\tilde g(s)\) the \(y_2\)-heights of the trace and the datum in the obstacle
chart, where \(s=y_1\) is the coordinate along the tangent line.
The rescalings \(u_r(z)=r^{-2}\bar u(rz)\)
converge to \(Q\) in \(C^1_{\mathrm{loc}}(\overline{\mathbb H_\alpha})\)
by~\cref{lem:boundary-blowups}. For \(\varepsilon>0\), we apply the barrier argument
in Step~2 of the proof of~\cref{prop:lipschitz} with \(v_\pm=\pm(1,a)+\varepsilon e_2\),
where \(a=\tilde g'(0)\). Indeed, for the constant \(K\) chosen there,
\(D_{v_\pm}Q=c\varepsilon n_2(z\cdot n)\ge KQ\) in a small half-ball, while
on the fixed line \(u_r=D_{(1,a)}u_r=0\), so that \(D_{v_\pm}u_r-Ku_r
=\varepsilon\partial_{z_2}u_r\ge0\). The comparison therefore gives
\(D_{v_\pm}u_r\ge0\) in a smaller half-ball for all small \(r\).
Scaling back, \(\{\bar u>0\}\) is the epigraph of \(s\mapsto as+\tilde h(s)-\tilde g(s)\),
so monotonicity in \(v_\pm\) implies that all tangent slopes of this graph
are between \(a-\varepsilon\) and \(a+\varepsilon\), giving
\(|\tilde h'_\pm(s)-\tilde g'(s)|\le\varepsilon\) for every \(s\) near zero.
Since \(\varepsilon\) is arbitrary and \(\tilde g'\) is continuous, we obtain
\(h'(p)=g'(p)\) and continuity of \(h'\) at \(p\).

To complete the proof of (ii) fix a point \(p\in\{h\neq g\}\setminus\Sigma\) and note that the inward curvature
is strictly negative by~\cref{prop:attachment}. Thus~\cref{lem:regular-points}
gives analyticity of the free boundary in the obstacle chart, and the
analytic change of coordinates yields the analiticity of $h$.
\end{proof}

\begin{remark}[Realization of a corner]\label{rem:corner}
The two branches in~\cref{thm:main}(iii) may meet at a corner. Indeed, let \(\Omega\) be a bounded
real-analytic domain whose boundary near a point \(p\) is concave, with
curvature \(\kappa(s)=-s^{2}+O(s^{3})\) in arclength (thus, \(\kappa\) vanishes to
second order at \(p\) and \(p\in\Sigma\)). Since the curvature is negative on both
sides of \(p\), the vertical-barrier construction of Simon~\cite{Simon1982} (see
also Williams~\cite{Williams1986}) furnishes analytic data \(g\) (for instance
a sufficiently large constant on the relevant arc) for which the datum is not
attained at \(p\), that is \(h(p)\neq g(p)\). The contact set then has positive
height over \(p\), and the free-boundary point above it is interior, with
forcing \(H(s)=c\,s^{2}+O(s^{3})\) and \(c>0\) (the
case \(m=1\)). By~\cref{lem:exclude-polynomial,lem:sector-classification} its
blow-up is the sector of opening \(\pi/3\), whose two rays make the angles
\(\pm\pi/6\) with the vertical.  Hence the free boundary has a corner at \(p\), pointing upward since \(g>h\) here.  Thus \(h\) is
right and left differentiable at \(p\) but not differentiable there, so the two
analytic branches need not have the same tangent.
\end{remark}

We now prove that if
\(p\in\Sigma\cap\{h\neq g\}\), then its two one-sided branches
extend real analytically across \(p\).
We use a variant of Kinderlehrer's Gauss-map reflection argument
for cylindrical obstacles~\cite[Section~7]{Kinderlehrer1973},
adapted to points where the curvature vanishes. More precisely, if the curvature
vanishes to order \(2m\), the boundary value of the Gauss map
vanishes to order \(2m+1\). We take this degeneracy into account by
introducing the coordinate \(w=\mathfrak g^{1/(2m+1)}\) and
applying reflection in the resulting sector.

\begin{proof}[Proof of~\cref{thm:main}(iii)]
We work in the coordinates of~\cref{sec:2d-blowups}, with the
free-boundary point at the origin.

\smallskip\noindent
\emph{Step 1: the Gauss map.}
By~\cref{prop:lipschitz}, we can choose
a small domain bounded by the two free-boundary arcs
near the origin, and a simple curve joining their other
endpoints through \(\{\bar u>0\}\).
Let \(Y=(Y_1,Y_2,Y_3)\) be a conformal
parametrization of the graph of \(\tilde u\) over this domain
by an upper half-disk, with the origin mapped to the origin
and the diameter mapped to the two free arcs.
The parametrization extends continuously to the diameter
by Caratheodory's theorem, and the normal is continuous
by~\cref{prop:c11}.

In these parameters, the Gauss map
\[
\mathfrak g
=\frac{-\partial_{y_1}\tilde u+i\partial_{y_2}\tilde u}
       {1+\sqrt{1+|\nabla\tilde u|^2}}
\]
is holomorphic and real on the diameter, while
\(\Im\mathfrak g>0\) in the upper half-disk
by~\cref{prop:strict-monotonicity}.
Hence Schwarz reflection extends \(\mathfrak g\) across the
diameter, and Hopf's lemma applied to \(\Im\mathfrak g\)
shows that its zero at the origin is simple.

\smallskip\noindent
\emph{Step 2: extension of the first coordinate.}
On the free boundary we have
\[
\mathfrak g=P(y_1):=
\frac{-\varphi'(y_1)}{1+\sqrt{1+\varphi'(y_1)^2}}
=c\,y_1^{2m+1}+O(y_1^{2m+2}),\qquad c>0.
\]
Let \(G\) be the analytic inverse of the signed real root
\(P^{1/(2m+1)}\), so that \(G(0)=0\) and \(G'(0)>0\).
We reparametrize the surface by
\[
w=\mathfrak g^{1/(2m+1)},\qquad
0<\arg w<\alpha:=\frac{\pi}{2m+1},
\]
and retain the notation \(Y_i\) for its coordinate functions.
Then
\[
Y_1(t)=G(t),\qquad
Y_1(te^{i\alpha})=G(-t)
\qquad\text{for all sufficiently small }t>0.
\]

We claim that \(Y_1\) extends harmonically to a neighborhood
of the origin. Indeed, write \(G(t)=\sum_{j\ge1}a_jt^j\)
and choose \(b_j\in\mathbb C\) satisfying
\[
\Re b_j=a_j,\qquad
\Re(b_je^{ij\alpha})=(-1)^j a_j.
\]
If \(j\) is a multiple of \(2m+1\), these equations are satisfied by \(b_j=a_j\),
since \(\cos(j\alpha)=(-1)^j\).
Otherwise they determine \(b_j\), and
\[
|\sin(j\alpha)|\ge\sin\alpha
\]
gives \(|b_j|\le C|a_j|\). Thus
\(B(w):=\sum_{j\ge1}b_jw^j\) is holomorphic near zero and
\(\Re B\) has the same boundary values as \(Y_1\) on both rays.
The bounded harmonic function \(Y_1-\Re B\) therefore extends
by successive odd reflections to a punctured disk.
The singularity at the origin is removable, which proves
the claim. In particular, \(Y_1=\Re F_1\) for a function
\(F_1\) holomorphic near zero.

\smallskip\noindent
\emph{Step 3: analyticity of the two branches.}
Writing \(Y_2=\Re F_2\), the Weierstrass relations give
\[
F_2'=-i\,\frac{1+w^{4m+2}}{1-w^{4m+2}}\,F_1',
\]
so that \(F_2\) also extends holomorphically through zero.
Hence the restrictions of \((Y_1,Y_2)\) to the two rays
are real analytic. Since their first coordinates are
\(G(t)\) and \(G(-t)\), with \(G'(0)>0\), the analytic
inverse function theorem shows that both free-boundary
branches extend as analytic graphs over \(y_1\).
The analytic change from \(y_1\) to arclength gives the
claimed extensions of the two branches of \(h\).
\end{proof}

\begin{remark}\label{r:Uraltseva}
The \(C^1\) conclusion of~\cref{thm:main}(ii) extends across a contact point
\(p\in\Sigma\) only if \(h=g\) near \(p\), as otherwise \(h\) has a corner there.
Indeed, if detachment accumulates at \(p\), one of the tangent rays of the trace
lies in the open half-plane by the blow-up argument in the proof of~\cref{lem:uniqueness-of-boundary-blow-ups}
and Step~3 of the proof of~\cref{prop:lipschitz}. Since \(h-g\) vanishes at \(p\) and has a fixed sign
nearby by~\cref{lem:vertical+separation}, differentiability would give
\(h'(p)=g'(p)\), forcing this ray to lie on the fixed line and contradicting
that it lies in the open half-plane.

Also, the \(C^1\) regularity near every point of nonzero
curvature on the relative boundary of \(\{h\ne g\}\) is expected to be optimal, since the classical obstacle
problem suggests that \(h'\) need not be Dini continuous at such points.
Indeed, consider the example described by Uraltseva
in \cite[Section 4]{Uraltseva1996}: it consists of a solution
\(w\in C^1(\overline{B_1^+})\) of
\[
\Delta w = \mathbf{1}_{\{w>0\}}
\quad\text{in } B_1^+,
\qquad
w(x_1,0)=0,
\]
where \(B_1^+:=\{(x_1,x_2)\in B_1:x_2>0\}\).
The free boundary \(\{w=0\}\cap B_1^+\) meets the fixed boundary
\(\{x_2=0\}\) at the origin, and although the closure of the free
boundary is a \(C^1\) curve, its normal fails to be Dini
continuous at the origin. Moreover, \(\partial_{x_2} w >0\) in
\(B_r\cap\{w>0\}\) for sufficiently small \(r>0\) (the \(C^1\) regularity of \(w\) implies that
\(\partial_{x_2} w=0\) on the closure of the free boundary, while
\cite[Theorem 2]{Uraltseva1996} gives
\(\partial^2_{x_2x_2} w \geq \frac12\) in \(B_r\cap\{w>0\}\)
for \(r>0\) small).

Nevertheless, one could still hope that \(\{h=g\}\) consists
of finitely many arcs, which are real analytic away from
finitely many exceptional points. This possibility is
consistent with Uraltseva's example, where the origin is
an isolated point of
\(\overline{\{w=0\}\cap B_1^+}\cap\{x_2=0\}\).
\end{remark}

\end{document}